\documentclass[a4paper,11pt,flqn]{article}
\usepackage{epsf}
\usepackage{enumitem}
\usepackage{graphicx}
\usepackage{mathptmx}
\usepackage{makeidx}
\usepackage{amsmath}
\usepackage{amssymb}
\usepackage{amsfonts}
\usepackage[bottom]{footmisc}
\usepackage{multicol}
\usepackage{type1cm}
\usepackage{mathrsfs}
\usepackage{url}
\usepackage{color}
\usepackage{lineno}
\usepackage[francais]{minitoc}
\usepackage{setspace}
\usepackage{float}
\usepackage{bm}
\usepackage{tikz}
\usepackage[hidelinks]{hyperref}
\usepackage[T1]{fontenc}
\usepackage[applemac]{inputenc}
\newtheorem{lemme}{Lemma \rm }
\newtheorem{prop}{Proposition \rm }
\newtheorem{rem}{Remark}

\newtheorem{definition}{Definition\rm}
\newtheorem{hyp}{Hypothesis\rm}
\newcommand{\fin}{\end{document}}
\newcommand{\beq}{\begin{equation}}
\newcommand{\feq}{\end{equation}}
\newcommand{\dcom}{\begin{quote}\begin{small}}
\newcommand{\fcom}{\end{small}\end{quote}}

\newcommand{\bq}{\begin{quote}}
\newcommand{\fq}{\end{quote}}

\newcommand{\Rmat}{\mathbb{R} }
\newcommand{\Nmat}{\mathbb{N} }

\newcommand{\Lmat}{\mathbb{L} }

\newcommand{\e}{\mathrm{e}}

\newcommand{\bitbul}{\begin{itemize}[label = \textbullet]}
\newcommand{\bittiret}{\begin{itemize}[label = -]}
\newcommand{\bito}{\begin{itemize}[label =$\circ$]}
\newcommand{\bit}{\begin{itemize}}
\newcommand{\fit}{\end{itemize}}
\newcommand{\ben}{\begin{enumerate}}
\newcommand{\fen}{\end{enumerate}}

\newcommand{\CN}{\mathcal{N}}
\newcommand{\chic}{\chi^{\mathcal{C}}}
\newcommand{\gr}{\textbf}
\newcommand{\dr}{\mathrm}
\newcommand{\implic}{\; \Longrightarrow \;}
\newcommand{\CNS}{\; \Longleftrightarrow \;}

\newcommand{\red}{\textcolor{red}}

\title{Sufficient condition for dispersal-induced growth on dynamic networks.}
\author{{{Michel Benaim}
\thanks{Institut de Math\'ematiques, Universit\' e de Neuch\^atel, Switzerland 
\tt{michel.benaim@unine.ch}}}
\and{{Claude Lobry} \thanks{C.R.H.I, Universit\'e de Nice Sophia Antipolis, France  \tt{claude\_lobry@orange.fr}  }}
\and{{Tewfik Sari}
\thanks{ITAP, Univ Montpellier, INRAE, Institut Agro, Montpellier, France
\tt{tewfik.sari@inrae.fr}}}
\and{{Edouard Strickler}
\thanks{Universit\' e de Lorraine, CNRS, Inria, IECL, Nancy, France$\quad \quad$
\tt{edouard.strickler@univ-lorraine.fr}}} }
\date{\today}
\begin{document}

\maketitle
\begin{tikzpicture}[remember picture,overlay]
\hypersetup{hidelinks}
\node[anchor=north west,yshift=280pt,xshift=-100pt]
{ \href{https://doi.org/10.24072/pci.mcb.100411}{\includegraphics[height=35mm]{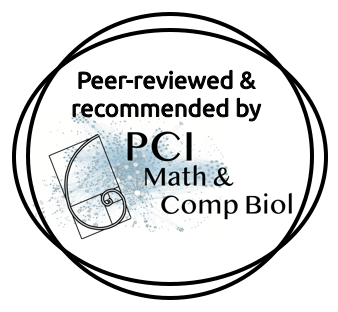}} } ;
\end{tikzpicture}
\begin{abstract}
We consider  a population spreading across a finite number of sites. Individuals can move from one site to the other according to a network (oriented links between the sites) that varies periodically over time. On each site, the population experiences a growth rate which is also periodically time-varying. Recently, this kind of models have been extensively studied, using various technical tools to derive precise necessary and sufficient conditions on the parameters of the system (ie the local growth rate on each site, the time period and the strength of migration between the sites) for the population to grow. In the present paper, we take a completely different approach: using elementary comparison results between linear systems, we give sufficient conditions for the growth of the population. This condition is easy to check and  can be applied in a broad class of examples. In particular, in the case when all sites are sinks (ie, in the absence of migration, the population becomes extinct in each site), we prove that when our condition of growth is satisfied, the population grows when the time period is large and for values of the migration strength that are exponentially small with respect to the time period, which provides a positive answer a conjecture stated by Katriel in  \cite{Katriel}.
\end{abstract}

%================
\section*{Introduction}
%================

An habitat where a population resides is called a ‘‘source'' when, in the absence of migration, environmental conditions ensure population growth, and called a ‘‘sink'' otherwise. When ‘‘sources'' and ‘‘sinks'' are connected in a network through which migrations occur, population growth on the network depends on various factors: the ‘‘environmental conditions'' at each site of the network, the structure of the network and the intensity of migrations. All these factors depend on time, in a more or less random way. Determining the growth conditions of a population on a network is a central theme in ecology, both for theoretical and practical reasons (see \cite{Baguette, Hanski}) described in \cite{Bennett} which we quote below :
\dcom
In the real world, patches of habitat vary greatly in the resources they provide for animals and in the disturbance they experience. Consequently, some populations can be regarded as ‘sources’ that produce a net surplus of animals that are available as potential colonists to other habitat patches. On the other hand, ‘sinks’ are those populations in which mortality exceeds natality and the persistence of the population depends on a regular influx of immigrants (...). There are, as yet, limited data on the relative frequency of sources and sinks in natural environments but theoretical models suggest that the relative proportions of each and the level of dispersal between them, may have a significant influence on regional population dynamics and species conservation.
\fcom
The present paper is a contribution to the theoretical models mentioned above. 

As early as 1997 and 1998 it was noticed by Holt \cite{HOLT97} and Jansen and Yoshimura \cite{Jansen} a paradoxical effect, later called {\em inflation} by Gonzalez and Holt \cite{ HOLT02} : \\\\
{\em When environmental conditions vary over time, it can happen that two habitats that are ‘‘sinks'' in isolation can become ‘‘sources'' when linked by migration.}\\\\
This paradoxical effect has motivated theoretical studies on continuous and discrete time models, deterministic and stochastic, for instance  \cite{KLA08,ROY05, SCH10, EVA13, SCH23} (see \cite{BLSSTPB},  \cite{BLSSJOMBa} and the references therein for more information).

Recently Katriel \cite{Katriel} and the authors of the present paper \cite{BLSSTPB, BLSSJOMBa,BLSSJOMBb,BLSSAFST} started a detailed mathematical study of the 
linear continuous time model when the local growth rate on each site is periodic \cite{Katriel}, periodic or stochastic   \cite{BLSSTPB, BLSSJOMBa,BLSSJOMBb,BLSSAFST}. In \cite{Katriel} Katriel  suggests renaming the inflation phenomenon {\em Dispersal Induced Growth} (DIG). Since this expression is more mathematically meaningful, we will  use it here.
  In all these articles, the inflation phenomenon is characterized in terms of the dominant eigenvalue $\lambda(m,T)$ of certain positive matrices associated with the model, depending on parameters such as the migration intensity $m$ and the period $T$. Thanks to mathematical results on certain symmetric positive operators \cite{Katriel}, Tychonov's theorem on singular perturbations of differential equations and Perron Frobenius' theorem \cite{BLSSTPB, BLSSJOMBa,BLSSJOMBb,BLSSAFST} it is then possible to describe the behavior of $\lambda(m,T)$ according to various  assumptions on the migration process (see also \cite{MSS24} for recent investigations). 

In the present paper, we take a completely different approach. Instead of trying to characterize the values of the parameters $m$ and $T$ for which inflation occurs, in a less ambitious way, we simply look for sufficient  conditions that can cause the phenomenon. 
This allows us to consider migration networks that are much more general, and therefore more realistic, than those imposed by the mathematical tools used previously. We consider a model where migration is described by a succession $\CN^1,\;\CN^2,\cdots$ of oriented graphs on a set of sites representing an evolution over the time of the structure of the migration network.  This approach, although absolutely elementary, explains, in our view, the reasons for the inflation phenomenon. Moreover it also answers a question raised in \cite{Katriel, BLSSTPB, BLSSJOMBa} about the order of magnitude of the $m$-threshold  at which the phenomenon appears.

 Organization of the article.\\
A suitable mathematical framework for describing our model is that of dynamic networks, increasingly considered in computer science and social network modelling (see for example \cite{CAST},\cite{DIB}) and recently appearing in theoretical ecology (see for instance \cite{KOL18}). Unfortunately, in the field of population dynamics, we know of no reference directly applicable to our situation. We therefore devote Section \ref{modele} to a detailed description of our model using notations close to those of graph theory. Section \ref{resultats} is the mathematical description of our central result,.
In Section \ref{applications} we apply our technique to solve  some questions raised in  \cite{Katriel, BLSSTPB,BLSSJOMBa,BLSSJOMBb} and, finally, in Section \ref{stochastique} we show briefly how our techniques apply to random systems.

Our results are direct consequence of Proposition \ref{minocircuit} which is the cornerstone of the present work.  This proposition is mathematically totally elementary, and we offer a detailed demonstration of it in an appendix which, for the reader's convenience, recalls some basic results on linear differential systems that can be found in any textbook.

\tableofcontents
\newpage
\section{Model and notations}\label{modele}

%=================
\subsection{The model}
%=================
%---------------------
\paragraph{Sites.}
We denote by
\beq
\Pi = \{ \pi_1, \pi_2, ...,\pi_i, ...,\pi_n \}
\feq 
a set of $n$ sites and by $x_i(t)$ the size of the population on the $i-$th site at time $t$.
We denote by 
$$\dr{x} = \big (x_1,..., x_i,..., x_n\big)$$
the vector of the $x_i\,'$s and conversely $[\dr{x}]_i$ is the $i-$th component of $\dr{x}$. The vector $\dr{x}$ is the meta-population. If $a = \pi_i$ is some site we denote $x_a(t) = x_i(t)$.
%------------------------------------------------
\paragraph{Equations of the dynamic.} 
We are interested in the system\\\\
\fbox{
\begin{minipage}{0.98\textwidth}
\beq \label{systemper1}
\Sigma(r_i(.) ,\ell_{ij}(.),m,T) = \left \{ \frac{dx_i}{dt} = r_i(t/T)x_i+ m \sum_{j=1}^n \ell_{ij}(t/T)x_j  \quad i = 1,...,n\right.\\[10pt]\quad
\feq
\end{minipage}
}\\\\
where 
\bito
\item the fonctions $r_i(t)$ and $\ell_{ij}(t)$  are  1-periodic, 
\item the matrix $L(t) = (\ell_{ij}(t))$ is a migration matrix which means that $\ell_{ij}(t) \geq 0 $ and  $\ell_{i,i}(t) = - \sum_{j \not = i}\ell_{ji}$ and thus $\sum_{\,i }\ell_{ij} = 0$. This last assumption is standing all over the paper and will not be repeated. We assume $m\geq 0$.
\fit
The system $\Sigma(r_i(.) ,\ell_{ij}(.),m,T)$ defined by \eqref{systemper1} is therefore a non-autonomous linear system of period $T$. The term $r_i(t/T)$ is the growth rate on the site  $\pi_i$ and the term $m\ell_{ij}(t/T)x_j $ is the migration rate from the site $\pi_j$ to the site $\pi_i$.\\\\
System $\Sigma(r_i(.) ,\ell_{ij}(.),m,T)$ is  a linear system of the form $\frac{d\dr{x} }{dt}= A(t)\dr{x}$ where $A(t)$ is
a matrix with positive off-diagonal elements (Metzler matrix).  Solutions of \eqref{systemper1} with positive initial conditions are positive (see Appendix \ref{metzler}).\\\\
\noindent We discuss the influence of parameters $T$ and $m$ on metapopulation growth. 
%---------------------------------
\paragraph{Link, network.}
A {\em link} on $\Pi$ is an arrow pointing from one site to another. It is noted
$$ \pi_i \to \pi_j$$
The link is {\em outgoing} from $ \pi_i$ and {\em incoming} in $\pi_j$. 

A set of links on $\Pi$ is a {\em network}, denoted $\CN$. In other words the pair  $(\Pi,\CN)$ is  a directed graph (di-graph) in the terminology of graph theory.
%------------------------------
\paragraph{Seasonality.}
The following assumptions are made.

The interval $[0,1]$ is the union of $p$ intervals
\beq
[0,1] = [t_0 =0, t_1]\cup[t_1,t_2]\cup \cdots\cup[t_{k-1},t_k] \cup\cdots \cup[t_{p-1},t_p = 1]
\feq

\begin{hyp}\label{H1}
All functions $r_i(t)$ and $\ell_{ij}(t)$ are \textbf{constant} on intervals $(t_{k-1},t_k) $ and we note
\beq
 r_i(t) = r_i^k\quad \quad \ell_{ij}(t) = \ell_{ij}^k\quad \mathrm{if} \quad t \in (t_{k-1},t_k) 
\feq  
\end{hyp}

The parameter $T$ is the duration of the periodic environment (e.g. year, day,...) and the successive intervals $[t_{k-1},\,t_k]$ are the successive ‘‘seasons''. Hypothesis \ref{H1} means that the growth functions $r_k(t/T)$ and migrations functions $\ell_{ij}(t/T)$ are constant during the $k-$th ‘‘season''. 

The upper subscripts in $r_i^k ,\ell_{ij}^k$ therefore indicate the season when these parameters are effective.
With these notations and an obvious interpretation in terms of ‘‘switched systems'' (see \cite{BLSSTPB}), we can rewrite the system \eqref{systemper1} as\\\\
 \fbox{
\begin{minipage}{0.98\textwidth}
\beq \label{systemper2}
\Sigma(r_i^k ,\ell_{ij}^k,m,T)=\quad \quad\left\{
\begin{array}{l}
\displaystyle \frac{dx_i}{dt} = r_i^kx_i+ m \sum_{j=1}^n \ell_{ij}^k x_j  \quad t \in [Tt_{k-1},\,Tt_k)\\[8pt]
i = 1,...,n\quad\quad k = 1,...,p 
\end{array}
\right.
\feq
$$\;$$
\end{minipage}}\\\\
which means that after having integrated the system \eqref{systemper2} up to time $Tt_k$ one takes $x_i(Tt _k);\; i = 1\cdot\cdot n$ as initial conditions for an integration on the interval $[Tt_k, \;Tt_{k+1}]$.
%--------------------------------------------------------------
\paragraph{Migrations on dynamic networks.}\label{migrations}

\begin{hyp}\label{zero-un}
For $i\not = j$,  $\ell^k_{ij} =0 \mbox{ or }1 $
\end{hyp} 

This hypothesis, which means that if there is migration between two sites its rate is always the same,  is rather restrictive.  Actually we make it in order to keep the things as simple as possible but it can be relaxed (see Subsection \ref{relax}).

\begin{figure}[h]
\center
\includegraphics[width=0.8\textwidth]{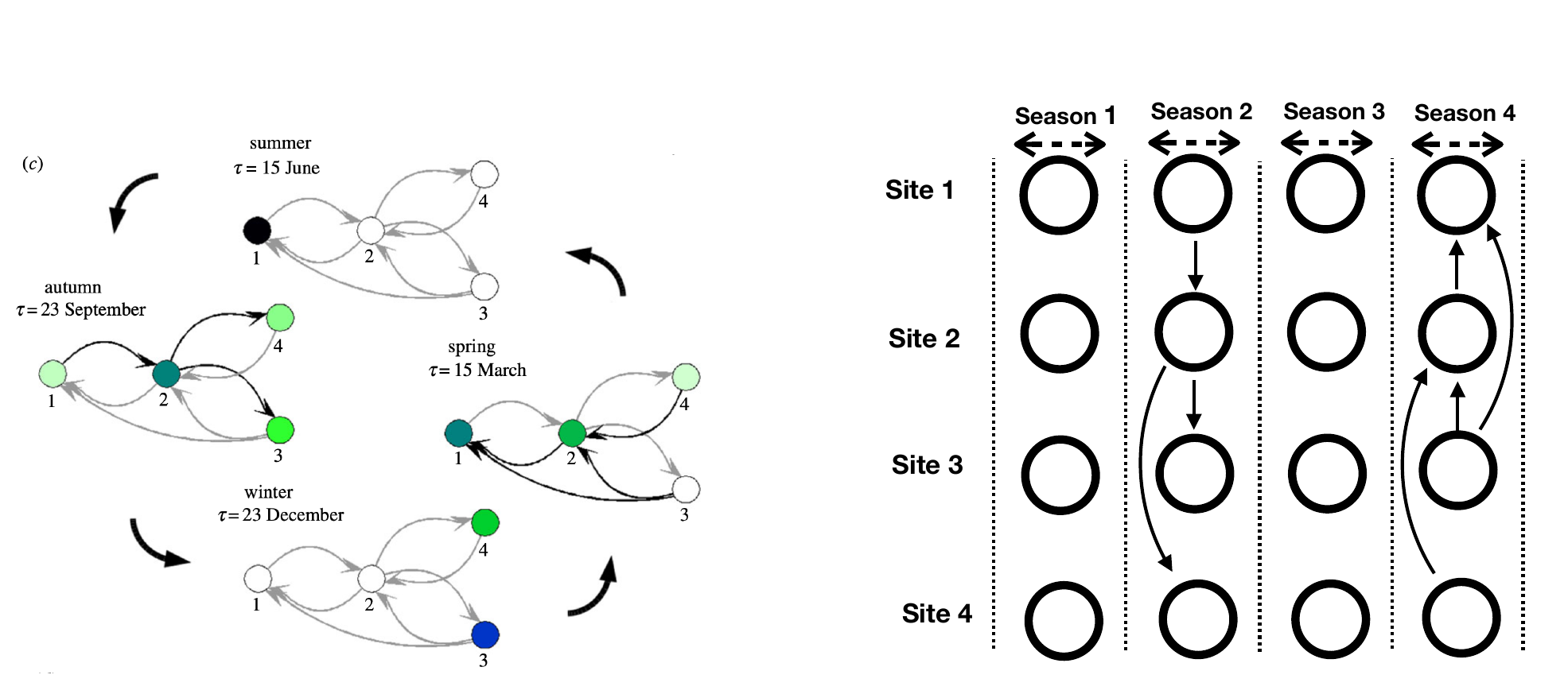}\\
\caption {A migration network from \cite{KOL18} .}\label{oiseaux}
\end{figure}

The $\left( \ell_{ij}^k \right) $ matrix is therefore equivalent to a  $\CN^k$ network on $\Pi$, defined by 
\beq
\pi_j \to \pi_i \,\in\, \CN^k \CNS \ell_{ij}^k = 1
\feq
The $\CN^k$ network is the {migration network} of the $k-$th season. The system \eqref{systemper2} is associated with the sequence 
\beq
\CN^{[1..p]} = \left\{\CN^1, ...,\CN^k,...,\CN^p \right \}
\feq
of $p$ networks.
\begin{definition}
An ordered set $\CN^{[1..p]} = \left\{\CN^1, ...,\CN^k,...,\CN^p \right \}$ of networks on the same set of sites  is called a $p$-dynamic network (or $p$-network in short). The migration network of the system $\Sigma(r_i^k ,\ell_{ij}^k,m,T)$ defined by the equations \eqref{systemper2} is its underlying p-network.
\end{definition} 
%--------------------------
\paragraph{Example.} 
In a paper entitled {\em  A periodic Markov model to formalize animal migration on a network}, A. Költz and al. \cite{KOL18} consider migration networks that change with seasons. In Figure \ref{oiseaux} we see on the left a scheme extracted from Figure 1 of the paper \cite{KOL18}. It represents four sites where birds are living.
In summer (season 1)  there is no migration (and most of the population is on site 1), in autumn there is migration to sites 3 and 4 through site 2, in winter (season 3) no migration and in spring (season four) migration back from sites 2, 3, 4 to site 1. On the right of our Figure 1 one sees the representation of this dynamic network in the style that we use in our paper.

%================================
\subsection{Terminology: paths, circuits.}\label{terminology}
%------------
\paragraph{Notations} We use the following notation: for $\dr{x},\,\dr{y}  \in \Rmat^N$, $\dr{x} ≥ \dr{y} $ means that for all $i$, $x_i \geq y_i $ ; $\dr{x} > \dr{y} $ means that $x_i \geq y_i $ and $\dr{x} \not = \dr{y} $ and $\dr{x} \gg \dr{y}  $ means that for all $i$, $ x_ i >y_i $.
 We use the same notation for $n \times p$ matrices considered as elements of $\Rmat^{n\times p}$.
\paragraph{Path} A  path (without loop) in a network $\CN$ is sequence of links $\gamma_i$, such that the origin of $\gamma_{i+1}$ is the extremity of $\gamma_i$, never returning to a previously visited site. Precisely:\\
 
 \noindent \fbox{
\begin{minipage}{0.98\textwidth}
A {\em path} $ a\,\Gamma\,b$ 
on the network $\CN$ is a sequence 
 \beq
a\Gamma b = \{a = \pi_{i(0)} \to\pi_{i(1)} \to \cdots\to \pi_{i(j)}\to \cdots \to \pi_{i(l)} = b\}
\feq
 of all \textbf{different sites} $\pi_{i(j)}$ of $\Pi$ linked by $\pi_{i(j-1)} \to \pi_{i(j)}$.  The index $l$ is the length  $\Lmat(a\Gamma b)$ of the path $ a\Gamma b$.
 \end{minipage}
}\\\\
If  there exists a path $a\Gamma b$ from the site $a$ to the site $b$, then $a$ is said to be ‘‘upstream'' of $b$ and $b$ is said to be ‘‘downstream'' of $a$.

Notice that what we call a {\em path} is a \textbf{path with no loop}. When we consider paths with loops later on, we will use {\em lpath} terminology.
%---------------------------------------
\paragraph{Dynamic  path.}
Let $\CN^{[1... p]} = \left\{\CN^1,...,\CN^k,  ...,\CN^p \right\}$ be a $p$-dynamic network. \\\\
\noindent \fbox{
\begin{minipage}{0.98\textwidth}
A  {\em dynamic $p$-path} of $\CN^{[1\cdots p]}$ is a sequence of paths
\beq
 a_0\Gamma^1a_1\Gamma^2a_2\cdots a_{k-1}\Gamma^ka_k\cdots a_{p-1}\Gamma^pa_p
\feq
where each path $a_{k-1}\Gamma^ka_k$ is a path of $\CN^k$. 
 \end{minipage}
}
\paragraph{Dynamic circuit.}
Let $\CN^{[1... p]} = \left\{\CN^1,...,\CN^k,  ...,\CN^p \right\}$ be a $p$-dynamic network. \\\\
\noindent \fbox{
\begin{minipage}{0.98\textwidth}

A dynamic $p$-circuit  ($p$-circuit in short) of the $p$-dynamic network $\CN^{[1..p]} $  is a dynamic $p$-path 
\beq
\label{eq:pcircuit}
 \mathcal{C} =  a_0\Gamma^1a_1\Gamma^2a_2\cdots a_{k-1}\Gamma^ka_k\cdots  a_{p-1}\Gamma^pa_p
\feq
such that $a_p = a_0$. The length of $\mathcal{C}$ is
$$ \Lmat(\mathcal{C} )= \sum_{i = 1} ^p \Lmat \left(a_{i-1}\Gammaîa_i \right)$$
 \end{minipage}
 }
 
\begin{figure}[h]
\center
\includegraphics[width=0.6\textwidth]{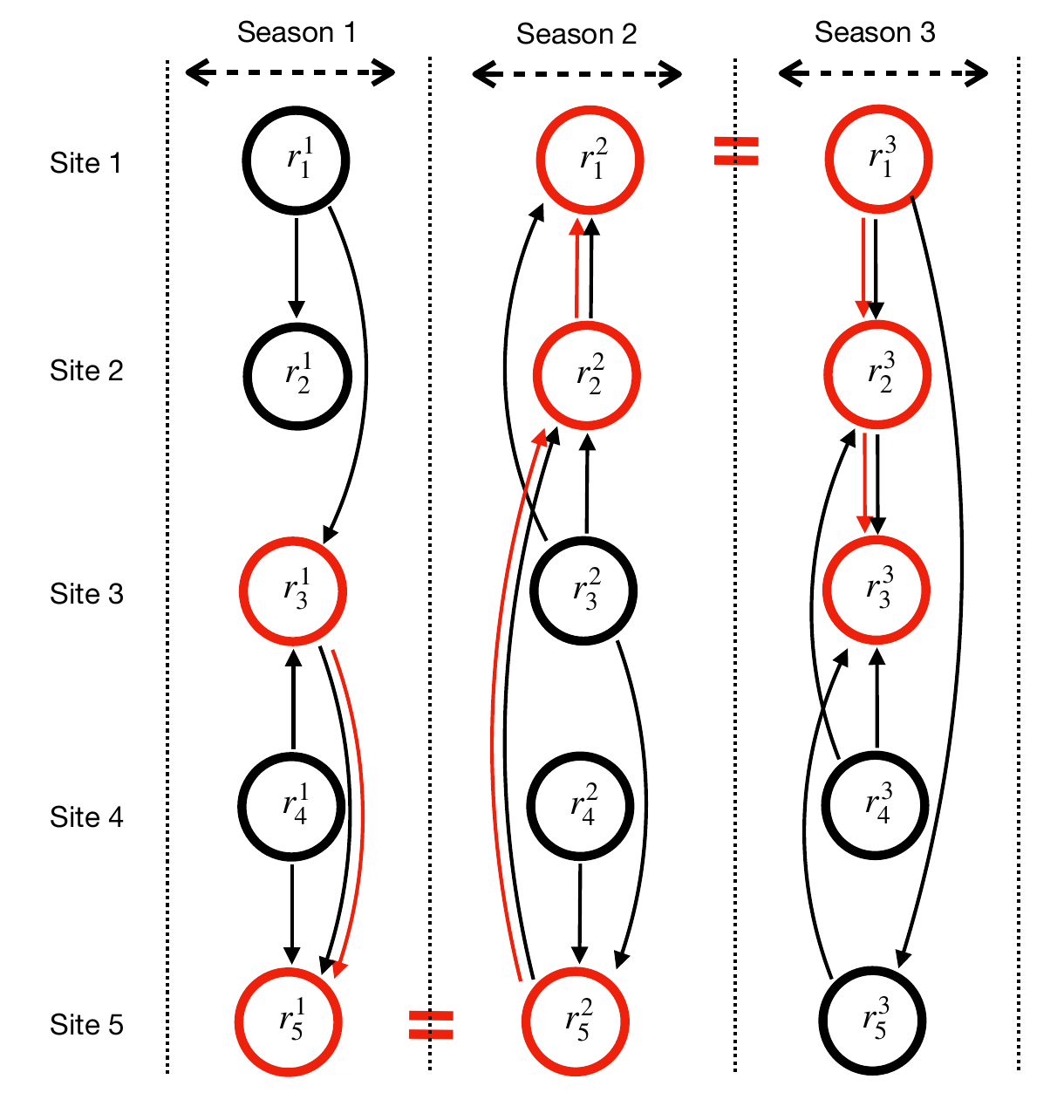}\\
\caption {In red a 3-circuit in a dynamic network (black arrows).}\label{reseau}
\end{figure}

 Figure \ref{reseau} shows an example of a dynamic  network (5 sites and 3 seasons) with a 3-circuit issued from site 3 in red.

%===============
\section{Main results}\label{resultats}
%===============
%=============================================
\subsection {Sources and sinks in a periodic system.}
Consider the $T$-periodic  switched system defined by \eqref{systemper2}
with  its $p$-underlying dynamic network
\beq \label{UN}
\CN^{[1..p]} = \left\{\CN^1, ...,\CN^k,...,\CN^p \right \}
\feq

\begin{definition}
One says that, in absence of migration, a site $\pi_i$ is a ‘‘source'' if $\sum_{k = 1}^p r_i^k(t_k-t_{k-1}) > 0$ which is equivalent to saying that the solution of the T-periodic switched system 
\beq \label{switch}
\frac{dx_i}{dt} = r^k_ix_i \quad \quad t \in [Tt_{k-1},\,Tt_k[ \quad \quad  x_i(0) = x_o > 0
\feq
tends to infinity as $t$ tends to infinity, since $x_i(t) = \exp( \int_0^t r_i(s)ds)x_o$. 

If $\sum_{k = 1}^p r_i^k(t_k-t_{k-1}) < 0$ the site $\pi_i$ is called a ‘‘sink'', which is equivalent to saying that the solution $x_i(t)$ goes to $0$ when $t$ goes to infinity.
\end{definition}
Notice that 
$$ \sum_{k = 1}^p r_i^k(t_k-t_{k-1}) = \frac{1}{T}\int_0^Tr_i(t) dt$$
is the average growth rate of the site $\pi_i$ in the absence of migration on the period.

We extend the definition of ‘‘source'' and ‘‘sink'' to the whole system $\Sigma(r_i^k ,\ell_{ij}^k,m,T)$  defined by \eqref{systemper2} by saying that it is a ‘‘source'' when, given an initial condition $\dr{x}(0) \gg 0$  the corresponding total population $S(t, m,T,\dr{x} (0))) = \sum_{i = 1}^n x_i(t,m,T,\dr{x}(0)) $ tends to infinity when $t$ tends to infinity. Since (due to linearity) the fact that $S(t,m,T,\dr{x}(0))$ tends to infinity is independent of the positive initial condition $\dr{x}(0)\gg 0$ we omit  it in the following.
\begin{definition}\label{seuil} \textbf{ $m^*(T)$-threshold.}  Let $S(t, m,T)$ be the total population of $\Sigma(r_i^k ,\ell_{ij}^k,m,T)$. The $m^*(T)$-threshold of $\Sigma(r_i^k ,\ell_{ij}^k,m,T)$ is the number 
\beq
m^*(T) = \inf \,\{m:\; S(t,m,T) \to \infty\}
\feq
\end{definition}

\subsection{A sufficient condition for Dispersal-Induced Growth (DIG)}
The DIG phenomenon refers to the fact that  the total population  growth rate can be higher than all the mean growth rates of the isolated sites. In particular we adopt the following definition from \cite{Katriel}.
\begin{definition}\textbf{{\em DIG \cite{Katriel}}}. We consider the system $\Sigma(r_i^k ,\ell_{ij}^k,m,T)$ defined by \eqref{systemper2} for which we assume that every site is a ‘‘sink'' (i.e. $\forall \,i:\sum_{k = 1}^p r_i^k(t_k-t_{k-1} )<  0)$. 
One says that there is DIG (Dispersal Induced Growth) for  $\Sigma(r_i^k ,\ell_{ij}^k,m,T)$  if there exist  $m,T$ such that  the whole system is a ‘‘source'' (i.e. $ S(t) \to \infty$) or, in other words, if there is some $T$ such that $m^*(T) < \infty$.
\end{definition}
In \cite{Katriel} Katriel introduced the ‘‘growth index'':
\begin{definition}\label{growth-index}
\textbf{Growth-index of the system.}\\The growth index of the system $\Sigma(r_i^k ,\ell_{ij}^k,m,T)$ defined by \eqref{systemper2} is the number
\beq
\chi = \sum_{k = 1}^p \left(\max_{i=1}^n r_i^k\right)(t_k-t_{k-1}) = \frac{1}{T} \int_0^T \max_{i = 1}^nr_i(t)dt
\feq
\end{definition}
Then Katriel proved the following 
\begin{prop}\label{CNDIG}
\textbf{\em - Necessary condition for DIG}  A necessary condition for the system $\Sigma(r_i^k ,\ell_{ij}^k,m,T)$ defined by \eqref{systemper2} to be a ‘‘source'' for some $m,T$, is that its growth index $\chi$ is positive. In other words, a necessary condition for the existence of DIG is that the growth index $\chi$ is positive.
\end{prop}
\textbf{Proof}. This is the proof given by Katriel in the preprint version \cite{KatrielV1}. The proof given in \cite{Katriel} relies on the monotonicity of the Lyapunov exponent as a function of $T$, which is not true in general when the migration matrix in not symmetric (see \cite[Proposition 3]{MSS24}). One has
\beq
\begin{array}{l}
\displaystyle  \frac{dS}{dt} = \sum_{i=1}^n \left( r_i(t) x_i(t) +m \sum_{j = 1}^n \ell_{ij}(t)x_j(t)\right)\\[8pt]
\displaystyle  \frac{dS}{dt} = \sum_{i=1}^n r_i(t) x_i(t) \leq \sum_{i=1}^n(\max_{i=1}^n r_i(t)) x_i(t) =(\max_{i=1}^n r_i(t))S(t)
  \end{array}
 \feq
Let $\rho(t) = \max_{i=1}^n r_i(t)$, one has $\frac{dS}{dt} \leq \rho(t) S(t)$ which implies $S(T) \leq \e^{\int_0^T \rho(s)ds}$ which, by definition of $\chi$ is  $S(T) \leq \e^{T\chi} S(0)$. Hence, if $\chi < 0$ the total population tends to $0$ and the system is not a ‘‘source''.$\Box$

\begin{rem}
Note that the condition $\chi > 0$ implies first that some of the habitats have, in some season, a positive growth rate. In addition, since each habitat is a sink, in order to have $\chi > 0$, the different habitats must be not too synchronised, in the sense that they shall not have all a positive growth rates in the same seasons and negative growth rates in the other seasons. For more details, see the discussion in Section 3.1 of \cite{BLSSTPB}.
\end{rem}

In \cite{BLSSJOMBa} we proved that whenever each network $\mathcal{N}^k$ is strongly connected\footnote{Meaning that for each pair of sites $a$ and $b$, there exists a path from $a$ to $b$} then the condition $\chi > 0$ is also sufficient for DIG. Nonetheless, if some of the networks are not strongly connected, we exhibited in~\cite{BLSSJOMBa} examples where $\chi > 0$ but DIG do not happen (see the example with 3 sites and 3 seasons in Section \ref{33} of the present paper). One of the main goal of the sequel is to give a sufficient condition for DIG\footnote{for other sufficient condition for DIG when some networks are not strongly connected, we refer to \cite{BLSSJOMBb}} relying on growth of the population on $p$-cicruit. To this end, we introduce now the {\em growth index of a p-circuit} which is the growth index of the system reduced on the $p$-circuit. Precisely:

\begin{definition} \textbf{Growth index of a $p$-circuit.} \\
Let 
\beq \label{circuit}
\mathcal{C} =
a_0\Gamma^1 a_1\Gamma^2 a_2 \cdots a_{k-1}\Gamma^k a^k \cdots
a_{p-1}\Gamma^pa_0
\feq
with 
$$ a_{k-1}\Gamma^ka_k = \left\{  a_{k-1} = \pi_{i^k(0)} \to \pi_{i^k(1)}\to \cdots \pi_{i^k(j)} \to \cdots \to \pi_{i^k(l^k)}= a_k\right \}$$
 be a  $p$-circuit defined on the underlying time-varying network of the system \eqref{systemper2}. We call {\em growth index} of the $p$-circuit $\mathcal{C}$ the number
\beq \label{growthindex}
\chi{\mathcal{C}} = \sum_{k = 1}^p (t_k - t_{k-1}) \max_{\pi_i \in \Gamma^k} r_i^k.
\feq
\end{definition}

\begin{rem}
Note that the growth index  $\chi^{\mathcal{C}}$ of a $p$-circuit is always smaller than $\chi$ since for all $k$,  $\displaystyle \max_{\pi_i \in \Gamma^k} r_i^k \leq \max_{i = 1\cdot \cdot n} r_i^k$. Moreover, when every network is strongly connected, we can always find a $p$-circuit $\mathcal{C}$ such that $\displaystyle \max_{\pi_i \in \Gamma^k} r_i^k = \max_{i = 1\cdot \cdot n} r_i^k$ for all $k$. Hence, for this circuit $\mathcal{C}$, $\chi^{\mathcal{C}} = \chi$ and the next proposition retrieves some of the results of the aforementioned paper \cite{BLSSJOMBa}.
\end{rem}

\begin{prop}\label{suffcroissance}
\textbf{{\em - Sufficient condition for DIG}}.
 Consider the T-periodic system $\Sigma(r_i^k ,\ell_{ij}^k,m,T)$ defined by \eqref{systemper2}.  Assume that there exists a  $p$-circuit  $\mathcal{C}$ of the underlying time-varying network of $\Sigma(r_i^k ,\ell_{ij}^k,m,T)$ such that $\chi^\mathcal{C} >0$. Then for all $\tau$ large enough, there exist   $0 < a(\tau) < b(\tau)$ such that that for every $T \geq \tau $, 
\begin{equation}
\label{eq:m-in-a-b}
 m\in[a(\tau),b(\tau)] \implic \;\Sigma(r_i^k ,\ell_{ij}^k,m,T)\mathrm{\; is \; a\;source}.
\end{equation}
In particular, the existence of a $p$-circuit $\mathcal{C}$ with a growth index $\chi^{\mathcal{C}} > 0$ is a sufficient condition for DIG.
\end{prop}
The condition in~\eqref{eq:m-in-a-b} can be rephrased as: {\em  If $\chi^\mathcal{C} >0$, then  for $T$ large enough and for $m$ neither too small nor too large, the system is a source.} 

The proof of Proposition~\ref{suffcroissance}  relies on the minoration of the growth given in the following proposition:
%----------------------------------
\begin{prop}\label{minocircuit}
Consider the system $\Sigma(r_i^k ,\ell_{ij}^k,m,T)$ defined by \eqref{systemper2}  on the underlying network $\mathcal{N}^{[1\cdots p]}$. 
Consider the $p-$circuit $\mathcal{C}$ defined by
\beq \label{circuit0}
\mathcal{C} =
a_0\Gamma^1 a_1\Gamma^2 a_2 \cdots a_{k-1}\Gamma^k a^k \cdots
a_{p-1}\Gamma^pa_0
\feq
with 
$$ a_{k-1}\Gamma^ka_k = \left\{  a_{k-1} = \pi_{i^k(0)} \to \pi_{i^k(1)}\to \cdots \pi_{i^k(j)} \to \cdots \to \pi_{i^k(l^k)}= a_k\right \}$$
and its growth index $\chi^{\mathcal{C}}$. 
 Then there exist $\tau$ and constants $C >0$ and $0 < \mu \leq n-1 $ (independent of $m$) such that 
\beq \label{minoration1}
T > \tau  \implic x_{i^1(0)}(T) \geq   C m^{\Lmat} \e^{T\left(\chi^{\mathcal{C}}- \mu \times m \right)} x_{i^1(0)}(0)
\feq
where $\Lmat = \Lmat(\mathcal{C})$ is the length of the circuit. 
\end{prop}
	
So, after one period, the value of the population size at the beginning site of the $p$-circuit is minorized by $C m^{\Lmat} \e^{T\left(\chi^{\mathcal{C}}- \mu \times  m \right)}$ times the initial value. If this number is greater than one then the size of the population is increasing.
The proof of this proposition, which is elementary but a bit intricate, is given in Appendix \ref{demprop}. From this proof, one sees that $\mu$ can be chosen equal to $n-1$, see Remark~\ref{rem:mu} just after the proof. But on some particular networks, we could take a $\mu$ strictly smaller than $n-1$. That is why we prefer to keep the constant $\mu$ indeterminate instead of $n-1$ in our statements.

We give here the idea behind the formal proof of Appendix \ref{demprop}.
\bito
\item We isolate the $p$-circuit $\mathcal{C}$ from the whole dynamic network by cutting  all the links incoming from outside of $\mathcal{C}$ and we add to each site of $\mathcal{C}$ outgoing links ‘‘to the clouds'' such that the total number of links leaving each site is $n-1$ ; each ‘‘link to the clouds'' can be considered as an  added mortality rate  $m$.
\item This new dynamic network defines a new system whose solutions minorize those of \eqref{systemper2}. This is straightforward since cutting incoming links suppresses positive terms in the right hand of \eqref{systemper2} and the addition of  out going links adds negative ones.
\item On each path $\Gamma^k$ of the $p$-circuit   the maximum $r^k$ of the growth rates  is attained in some dominant site $\pi_{i^k(d)}$ and it is possible (Lemma \ref{lemme} of Appendix \ref{demprop}) to
prove by direct computation that for large enough $T$, for all the sites which are down stream of  $\pi_{i^k(d)}$ the growth rate is $r^k$. Actually this is straightforward  in the case of two sites. In this case the system is:

\begin{minipage}{0.45\textwidth}
$
\begin{array}{lcl}
\displaystyle \frac{dx_1}{dt} &=&\displaystyle( r_1-m)x_1\\[8pt]
\displaystyle \frac{dx_2}{dt} &=&\displaystyle mx_1 +(r_2-m)x_2
\end{array}
$
\end{minipage}\quad \quad
\begin{minipage}{0.15\textwidth}
\includegraphics[width=1\textwidth]{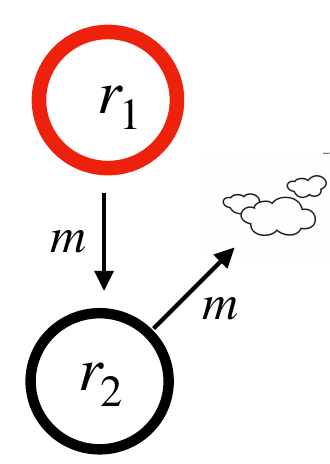}\\
\end{minipage}\\
Integration of the first equation gives
$ x_1(T) = \e^{T( r_1-m)}x_1(0) $
and the second one (see Appendix \ref{eqlineaire}) 
$$
\begin{array}{l}
x_2(T) = \e^{T(r_2-m)}\left\{ x_2(0) + \int_0^T \e^{-s(r_2-m)} m  \e^{s( r_1-m)}x_1(0) ds\right\}\\[8pt]
x_2(T)\geq  
 \e^{T(r_2-m)}\int_0^T \e^{s(r_1-r_2)}ds\, m x_1(0)\\[8pt]
x_2(T) \geq  \e^{T(r_2-m)}\frac{1}{r_1-r_2}\left[ \e^{T(r_1-r_2)}-1\right]\, m x_1(0) = \e^{T(r_1-m)}\frac{1-\e^{T(r_2-r_1)}}{r_1-r_2}\, mx_1(0).
\end{array}
$$
Since $r_1 > r_2$ the term $\e^{T(r_2-r_1)}$ tends to $0$ when $T$ tends to $\infty$ and therefore is smaller than $1/2$ for $T$ large enough and thus
$$ \exists \tau \; \mathrm{such\; that} \; T \geq \tau  \Longrightarrow x_2(T) > \frac{1}{2(r_1-r_2)} m \e^{T(r_1-m)} x_1(0)$$
which looks like \eqref{minoration1} with $\mu = 1$.
Proving the general case is just a matter of notations.
\item The iteration of this inequality along the paths of the $p$-circuit gives rise to \eqref{minoration1}.
\fit
\subsection{Relaxation of hypotheses}\label{relax}
%-------------------
\subsubsection*{Relaxation of Hypothesis 1}\label{relax1}

From a modeling perspective, Hypothesis 1 means that environmental parameters remain fixed throughout an entire season, which is obviously somewhat restrictive. However, by subdividing a season into multiple “mini-seasons,” we can more closely approximate the real-world situation.

From a mathematical point of view, a continuous function is the limit of piecewise constant functions, and there is no doubt that the “piecewise constant” assumption could be replaced by “piecewise continuous” in most of our results. This work remains to be done.

\subsubsection*{Relaxation of Hypothesis 2}\label{relax2}
Consider the system 
\beq \label{systemper3}
\Sigma(r_i^k,\alpha_i^k,\ell_{ij}^k,m,T)=\quad \quad\left\{
\begin{array}{l}
\displaystyle \frac{dx_i}{dt} = (r_i^k- m\alpha_i^k)x_i + m \sum_{j=1}^n \ell_{ij}^k x_j  \quad t \in [Tt_{k-1},\,Tt_k[ \\[8pt]
i = 1,...,n\quad\quad k = 1,...,p \quad \quad \alpha_i^k \geq 0
\end{array}
\right.
\feq
where the growth rate on each site decreases linearly with the parameter $m$.  It is clear that the system 
\beq \label{systemper4}
\Sigma(r_i^k,\alpha,\ell_{ij}^k,m,T)=\quad \quad\left\{
\begin{array}{l}
\displaystyle \frac{dx_i}{dt} = (r_i^k- m\alpha)x_i+ m \sum_{j=1}^n \ell_{ij}^k x_j  \quad t \in [Tt_{k-1},\,Tt_k[ \\[8pt]
i = 1,...,n\quad\quad k = 1,...,p \quad \quad \alpha \geq 0
\end{array}
\right.
\feq
with $\alpha = \max_{i,k}\alpha_i^k$, minorizes the system $\Sigma(r_i^k,\alpha_i^k,\ell_{ij}^k,m,T)$. So Proposition \ref{minocircuit} is still true for this system with a new  $\mu$ equal to $\mu +\alpha$.

Now we do not assume that $\ell_{ij} \in\{0,1\}$  in the system   $\Sigma(r_i^k ,\ell_{ij}^k,m,T)$ defined by \eqref{systemper2}. Let 
$$\ell = \min_{\ell^k_{ij} >0} \ell^k_{ij}; \; \quad \quad \quad \;   \beta^k_{i,j} = 1\; \mathrm{if} \; \ell^k_{ij} > 0, \; \beta^k_{i,j} = 0\;\mathrm{otherwise.}$$
 For $i \not = j$ we have $m \ell^k_{ij} \geq m \ell \beta^k_{i,j} $ and the system
\beq \label{systemper3}
\Sigma(r_i^k,\alpha_i^k,\beta_{i,j}^k,m,T)=\quad \quad\left\{
\begin{array}{l}
\displaystyle \frac{dx_i}{dt} = (r_i^k- m\alpha_i^k)x_i + m \ell \sum_{j=1}^n \beta_{ij}^k x_j  \quad t \in [Tt_{k-1},\,Tt_k[ \\[8pt]
i = 1,...,n\quad\quad k = 1,...,p 
\end{array}
\right.
\feq
with $\beta^k_{i,i} = -\sum_{j\not = i} \beta_{j,i}^k$ and $\alpha_i^k = \sum_{j\not = i} (\ell_{j,i}^k - \ell \beta_{j,i}^k) \geq 0$, to which we can apply Proposition \ref{minocircuit},  minorizes  the system $\Sigma(r_i^k ,\ell_{ij}^k,m,T)$. Thus we can relax hypothesis \ref{zero-un}.

%====================================
\subsection{The shape of the minimizing function.}\label{minimizingfunction}

Let us see now  how Proposition \ref{suffcroissance} follows immediately from Proposition \ref{minocircuit}. Indeed from \eqref{minoration1}, there exist $\tau$, $C$ and $\mu$  such that 
$$T > \tau \implic x_{i^1(0)}(T) \geq C m^{\Lmat} \e^{T(\chi^{\mathcal{C}}- \mu\times m)}x_{i^1(0)}(0)$$
Thus
 $C m^{\Lmat} \e^{T(\chi^{\mathcal{C}}- \mu \times m)} > 1$
implies that the sequence 
$ j \in \Nmat  \mapsto x_{i^1_0}(j\times T)$ 
tends to infinity which means that 
$\Sigma(r_i^k ,\ell_{ij}^k,m,T)$ 
is a ‘‘source''.

We let
$$H(m,T)\stackrel{\mathrm{def}}{=}   m^{\Lmat} \e^{T(\chi^{\mathcal{C}}- \mu \times m)} $$
denote the minimizing function appearing in \eqref{minoration1}. From elementary calculus it follows that when $\chi^{\mathcal{C}} > 0$,  \\
\begin{minipage}{0.45\textwidth}
\includegraphics[width=1\textwidth]{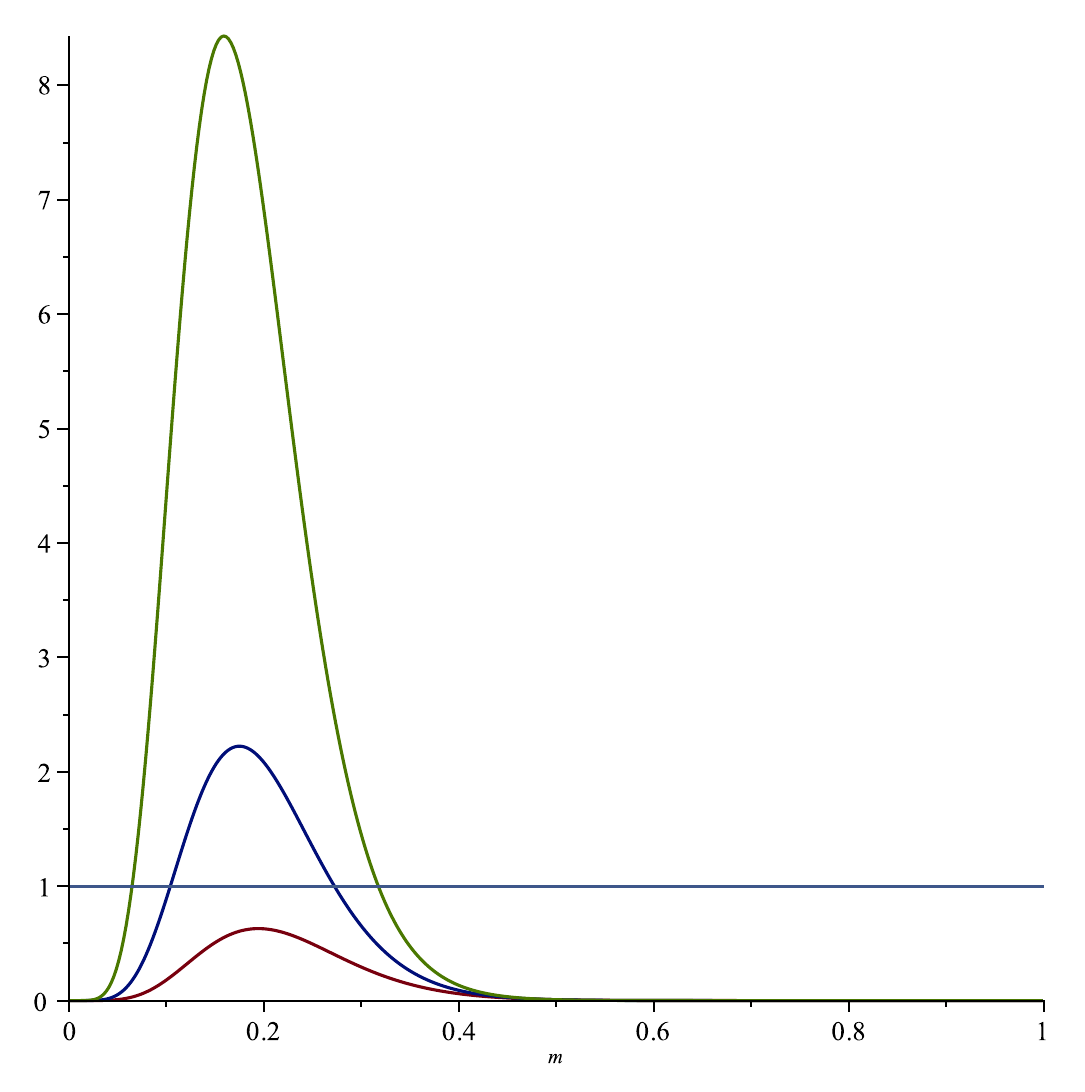}\\
\end{minipage}$\quad$
\begin{minipage}{0.5\textwidth}
\bito
\item   $H(m,T)$ is positive
\item  For $m < \frac{\chi^{\mathcal{C}}}{\mu}$, the function $T \mapsto H(m,T)$ is increasing,
\item  the function $m \mapsto H(m,T)$ is increasing, passes through a maximum $H_{\max}(T)$ and then decreases to~$0$,
\item  moreover $H_{\max}(T) \to \infty$ when $T \to \infty$.
\fit 
\end{minipage}\\
(Above the graphs  of the functions $m \mapsto H(m,T) =  m^7\e^{T(2-4m)}$ for $T = 9$ (red), $T= 10$ (blue), $T= 11$ (green).)

Thus, for $T$ large enough,  there is an interval $[a(T), b(T)]$  where $H(m,T) > 1$. which proves  Proposition \ref{suffcroissance}. 
%========================
\subsection{An exemple.}
\begin{minipage}{0.45\textwidth}
 \includegraphics[width=0.9\textwidth]{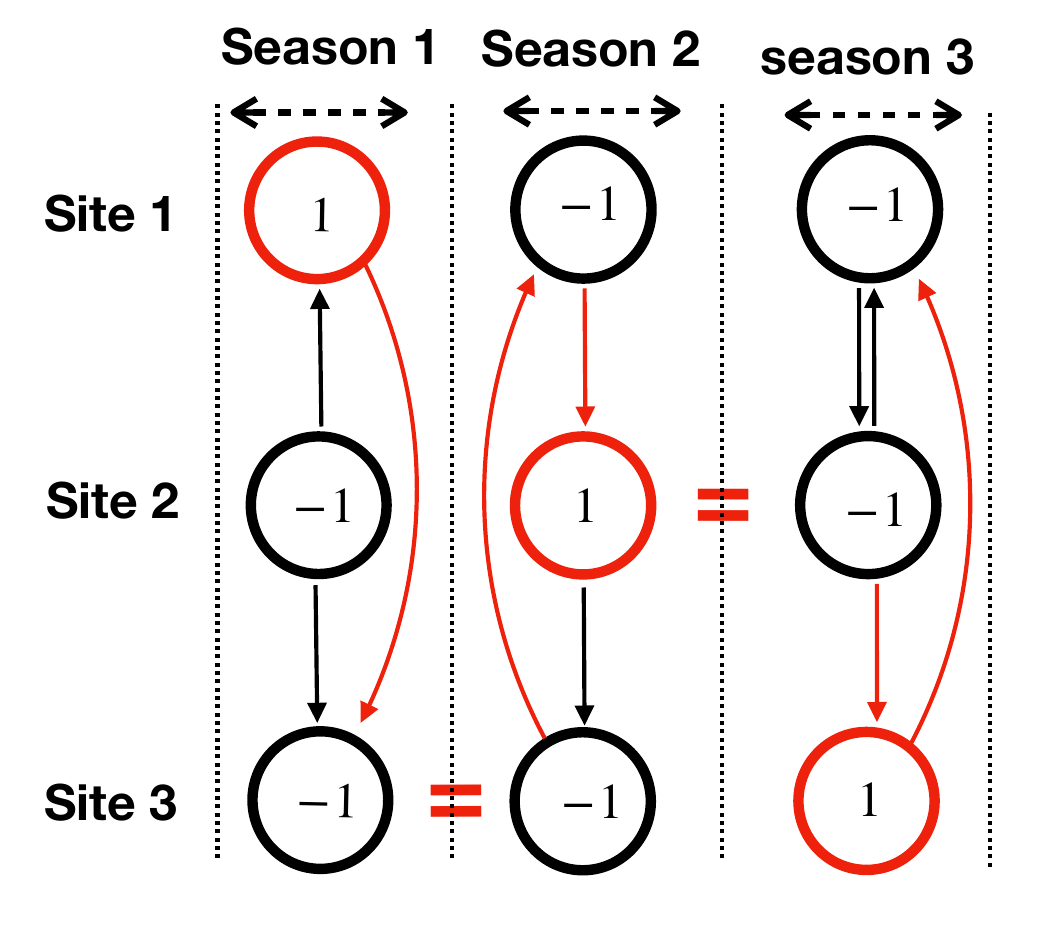}
\end{minipage}$\quad$
\begin{minipage}{0.5\textwidth}
Consider the system $\Sigma $ defined by the scheme on the left and assume that each season is of duration $1/3$. Each isolated site is a ‘‘sink'' (one  sequence with growth rate equal to 1 against two with a decay rate of $-1$). Consider the $3$-circuit $\mathcal{C}$  shown in red, that is:\\
$\mathcal{C} =\stackrel{\mathrm{season 1}}{|1\to 3|}\stackrel{\mathrm{season 2}}{|3\to 1\to 2|}\stackrel{\mathrm{season 3}}{|2\to 3 \to 1|}$
\end{minipage}\\
Since in each season the growth rate of the dominant site is $+1$ the  growth index of the circuit is $\frac{1}{3}+\frac{1}{3}+\frac{1}{3} = 1$ and thus is strictly positive. From Proposition \ref{suffcroissance} there is possibility of DIG.
\begin{figure}[h]
\center
\includegraphics[width=1\textwidth]{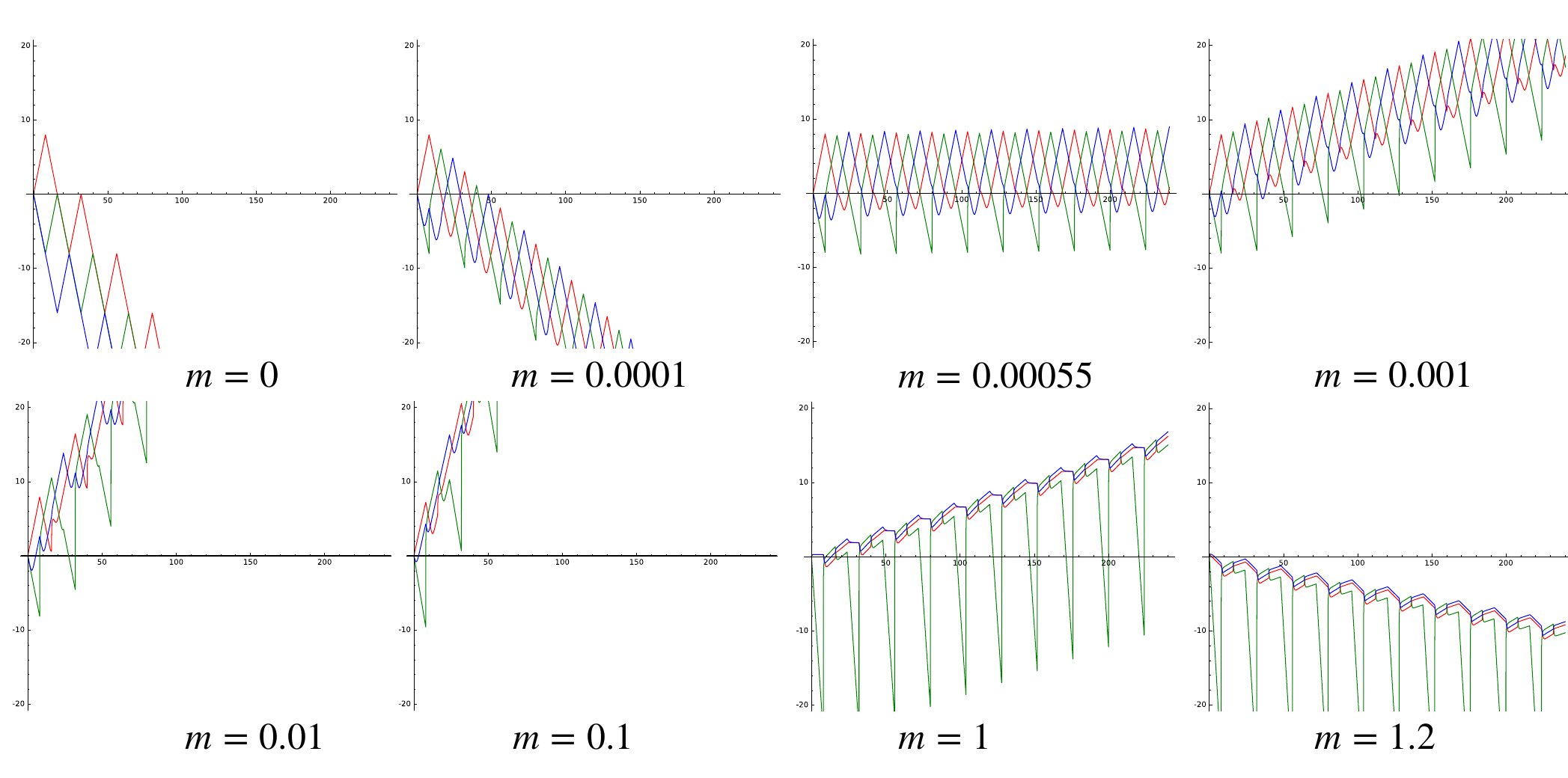}\\
\caption {Simulations of $\Sigma$ in the case $T =24$. One sees the logarithm of the population on site 1 (in red), on site 2 (in green) and site 3 (in blue)}\label{ex1}
\end{figure}
We have simulated the system $\Sigma$ from the initial condition $(1,1,1)$ in the case $T = 24$ and increasing values of $m$. On Figure  \ref{ex1} one sees the logarithm of  the population size on each site as a function of time. 
\bito
\item For $m=0$ all the sites are disconnected.  Unsurprisingly, since we have taken the logarithm, growth at site 1 is a straight line with a slope of +1 for 1/3 of the period, followed by a slope of -1 for the remaining two-thirds. The same occurs with the two other sites with a shift in the growing period.
\item For $m = 0.0001$ the total population is still decaying, while for $m = 0.001$ it is growing with a threshold around $m = 0.00055$.
\item The rate of growth increases for $m= 0.01$  and $m = 0.1$  and decreases later for $m = 1$
\item For $m = 1.2$ one sees that the system is again a ‘‘sink''.
\fit
An interesting point in this example is that the first value for which one observes an increase of the population is around $0.00055$ which is rather small compared to the parameters. This is actually a general feature that we explain in the next paragraph.
%========================================================
\subsection{The m-threshold is exponentially small with respect to $T$.}

Recall that we defined (see Definition \ref{seuil}) the m-threshold as the first value of $m$ such that, for a given $T$, the system $\Sigma(r_i^k ,\ell_{ij}^k,m,T)$ defined by \eqref{systemper2} is a source.
$$
m^*(T) = \inf_m \{ m\; \mathrm{such \;that \;} \Sigma(r_i^k ,\ell_{ij}^k,m,T)\;\mathrm{is\;a\;source} \}
$$
By Proposition~\ref{suffcroissance}, for $T$ large enough,  $m^*(T) \leq a(T)$  and in particular, $m^*(T) < \infty$.
\begin{prop}\label{expetit}
 Consider the T-periodic system $\Sigma(r_i^k ,\ell_{ij}^k,m,T)$.  Assume that there exists a  p-circuit  $\mathcal{C}$ of the underlying time-varying network of  $\Sigma$ such that $\chi^{\mathcal{C} }>0$. Then there exist $\tau$ and $\alpha > 0 $  such that  
 $$ T >\tau \implic \; m^*(T) \leq   \e^{ -\alpha T}$$
\end{prop}
\textbf{Proof}. From Proposition \ref{suffcroissance} we know that $\Sigma(r_i^k ,\ell_{ij}^k,m,T)$ is growing provided
$$ H(m,T) = C m^{\Lmat} \e^{T(\chi^{\mathcal{C}}- \mu m)} > 1$$ When $m < \frac{\chi^{\mathrm{C}} }{2\mu}$ we have $H(m,T) > C m^{\Lmat} \e^{T\frac{ \chi^{\mathcal{C}}}{2}  }$ and thus 
$$C m^{\Lmat} \e^{T\frac{ \chic}{2}  }    > 1 \implic   \Sigma(r_i^k ,\ell_{ij}^k,m,T)\;          \mathrm{is\; growing}$$
One easily checks that
$$ m\geq \e^{-T \frac{ \chic}{ 2\Lmat C^{1/\Lmat}}  }   \implic \; C m^{\Lmat} \e^{  T\frac{ \chic }{2}  } > 1$$
which proves the proposition with $\alpha = \frac{ \chic}{ 2\Lmat C^{1/\Lmat}}$.$\Box$\\

This proposition says that the $m$ threshold for DIG is exponentially small with respect to $T$. It  answers positively to the question asked in \cite{Katriel, BLSSJOMBa}.
%=================== 
\section{Some applications}\label{applications}
\subsection{The case of irreducible migration matrices}

In \cite{BLSSJOMBa} we considered the system 
$$
\Sigma(r_i(.) ,\ell_{ij}(.),m,T)\quad \quad \quad \quad \frac{dx_i}{dt} = r_i(t/T)x_i+ m \sum_{j=1}^n \ell_{ij}(t/T)x_j  \quad i = 1,...,n\\[10pt]
$$
where  $r_i(\tau) $ and $\ell_{ij}(\tau)$ are piecewise continuous functions of period 1, in the case where the migration matrix $L(\tau) =\left(\ell_{ij}(\tau)\right) $ is irreducible for every value of $\tau$ (this assumption, means that at each time, every patch is reachable from every other patch, either directly or by a path through other patches). 
	We proved (as a consequence of Perron-Frobenius theorem) that the Lyapunov exponent 
	$$\Lambda(m,T)[x_i] \stackrel{\mathrm{def}}{=} \lim_{t \to \infty} \frac{1}{t}\log(x_i(t))$$
	is actually independent of the site $i$ and, consequently, is denoted $\Lambda(m,T)$. The description of the growth of the system in terms of $\Lambda(m,T)$ is more precise than the one we propose here by just minorizing the growth but relies on the stronger assumption that $M(\tau)$ is irreducible. Nevertheless, our Proposition \ref{expetit} answers positively (at least in the piecewise constant case) to a question that was already asked in  \cite{Katriel, BLSSJOMBa}: {\em Assume that the growth index $\chi$ is positive. Is the growth threshold exponentially small with respect to $T$ ?} 
	Indeed, if for any $t$ the migration matrix is irreducible, during any season any site is connected to every site by a patch.  Then, starting from a dominant site at time $0$ there exists a  path during season 1 that connects it to one of the dominant sites of season 2, then, during season 2 there exists a  path that connects it to a dominant site of season 3, and so on until the last season when we return to the starting site; this defines a $p$-circuit whose growth index is precisely $\chi$ which allows us to apply the  Proposition \ref{expetit} which gives a positive answer to the conjecture.
	
	\begin{rem}
	\label{rem:totalmigration-Lambda}
In \cite{BLSSAFST, BLSSJOMBb}	we prove that the top Lyapunov exponent $\Lambda(m,T)$ is well defined and independent of $i$ as soon as the \textbf{total migration matrix} $L := \int_0^1 L(\tau) d \tau$ is irreducible. This assumption means that if we consider the network $\mathcal{N}^{tot}$ such that 
\[
\pi_j \to \pi_i \in \mathcal{N}^{tot}  \CNS \int_0^1 \ell_{ij}(\tau) d \tau > 0,
\]
then for every sites $a$ and $b$, there exists a path in $\mathcal{N}^{tot}$ from $a$ to $b$. In other words, the pair $(\Pi, \mathcal{N}^{tot})$ defines a strongly connected directed graph in the terminology of graph theory. Note that Proposition~\ref{expetit} applies in that situation, and thus, the assumption that all the migration matrices are irreducible is not needed to conclude that in case when there exists a $p$-circuit with positive growth index, the growth threshold is exponentially small with respect to $T$.
	\end{rem}
%=======================
%----------------------------------------------------
\subsection{Migratory birds example.}
This example is one of the simplest examples of growth on a time-varying network that one can imagine.

Consider for instance seasonally migrating species 
such as storks. It can be idealized  by the following model.
\begin{figure}[h]
\center
\includegraphics[width=0.8\textwidth]{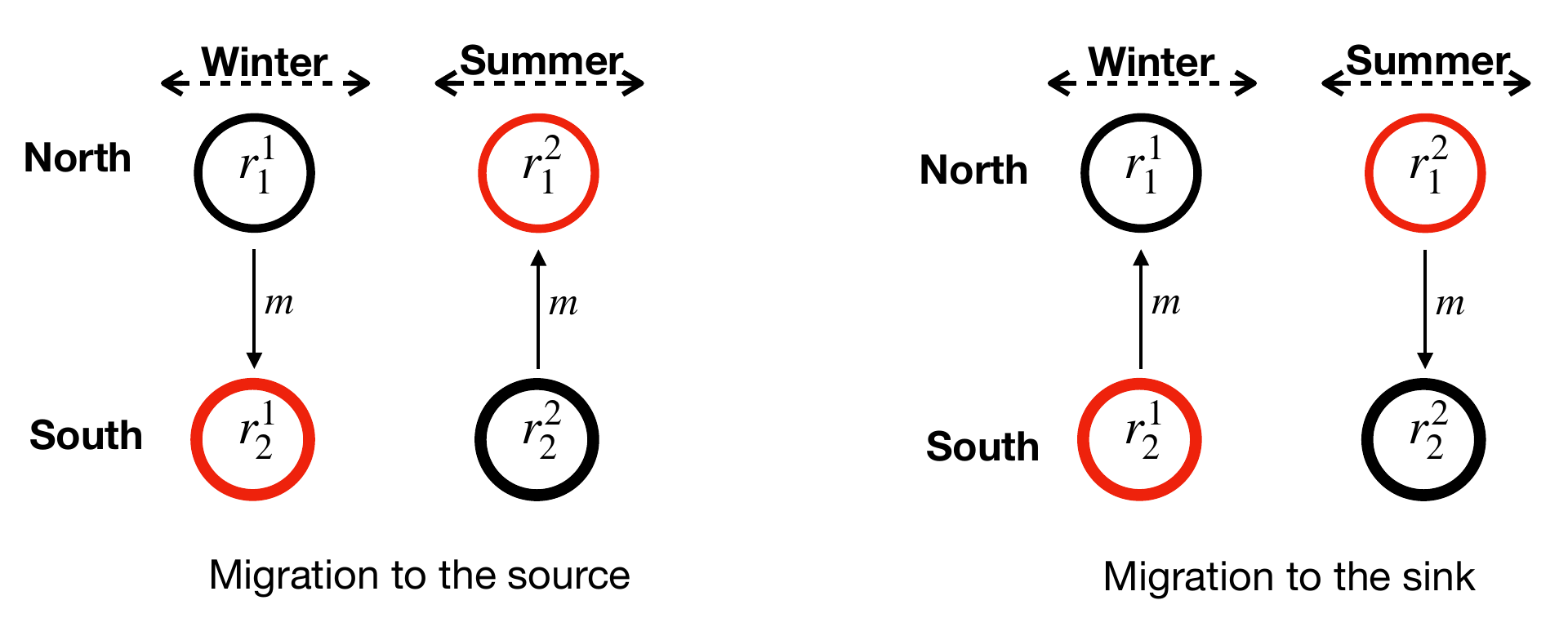}\\
\caption {Two possibilities for migration.}\label{cig}
\end{figure}
We consider two sites, say $\pi_1$ for north, $\pi_2$ for south and  consider two seasons, say winter and summer. Assume that in winter the south is a ‘‘source'' and the north is a  ‘‘sink'' while in summer the opposite is the case. 
%==========
\subsubsection*{Migration to the source}On the Figure \ref{cig}, on the left, we indicate a migration from north to south in winter and vice versa in summer which could be called a ‘‘migration to the source''. The corresponding equations are
\beq\label{vlS}
\mathrm{Winter} \left\{ 
\begin{array}{lcr}
\displaystyle \frac{dx_1}{dt} & = & (r_1^1-m)x_1\\[6pt]
\displaystyle \frac{dx_2}{dt} & = & mx_1+ r_2^1x_2
\end{array}
\right.
\quad
\mathrm{Summer} \left\{ 
\begin{array}{lcr}
\displaystyle \frac{dx_1}{dt} & = & r_1^2x_1+mx_2\\[6pt]
\displaystyle \frac{dx_2}{dt} & = & ( r_2^2-m) x_2
\end{array}
\right.
\feq
where $r_1^1 < 0 < r^1_2 $ and $r_2^2 <  0 < r_1^2 $.\\
$\mbox{‘‘North'' } \to \mbox{ ‘‘south'' } \to \mbox{ ‘‘north'' } =\mathcal{C}_{nsn}$ is a 2-circuit such that $\chi^{\mathcal{C}_{nsn}} > 0$. From  Proposition \ref{suffcroissance} we expect that the whole system is a source with $T$ large enough and  not too small not too large, $m$. 

 But for this  very simple system is is possible (using formal calculus software) to compute explicitly the growth. Indeed, let
$r_1^1 = r_2^2 = s, \,r_2^1 = r_1^2 = 1 $ with $s<0$ , and assume that the two seasons are of equal length $\frac{T}{2}$ and let
\beq
A^1 = 
\left(
\begin{array}{cc}
s-m&0\\
m&1
\end{array}
\right)
\quad \quad
A2= \left(
\begin{array}{cc}
1&m\\
0&s-m
\end{array}
\right)
\feq
 the matrix of the linear system operating during winter and during summer respectively. After winter the population is given by 
$$\dr{x}(T/2) = \e^{\frac{T}{2} \cdot A^1}\dr{x}(0)$$
and after summer by
$$\dr{x}(T) = \e^{\frac{T}{2}\cdot A^2}\dr{x}(T/2)=  \e^{\frac{T}{2}\cdot A^2} \e^{\frac{T}{2} \cdot A^1}\dr{x}(0)$$
The matrix 
$$M(m,T) \stackrel{\mathrm{def}}{=} \e^{\frac{T}{2}\cdot A^2} \e^{\frac{T}{2} \cdot A^1}$$
(called the monodromy matrix of the periodic system) expresses the growth after one period.  
Let  $\lambda(m,T)$ be  its dominant eigenvalue (i.e. with the greatest modulus). One knows that the system \eqref{vlS} is a source if $ \lambda(m,T)> 1$  and a sink if $ \lambda(m,T)< 1$. Moreover one knows that $\frac{\log(\lambda(m,T))}{T}$ is the asymptotic growth rate (i.e. the Lyapunov exponent) of $x_1(t) $ and $x_2(t) $:
$$\frac{\log(\lambda(m,T))}{T} = \lim_{t \to \infty}\frac{\log(x_1(t))}{t} = \lim_{t \to \infty}\frac{\log(x_2(t))}{t}  $$
In principle it is not difficult to compute $M(m,T)$ and its eigenvalues but it is a bit intricate and it is better to ask to some software (here we used Maple)  to compute them for us.
Maple says that the dominant eigenvalue is given by 
$$\lambda(m,T) = \frac{ \left(A+\sqrt{(B+C+D+E)m^2  }\right)}{2 \left(m -s +1\right)^{2}}$$
$$
\begin{array}{l}
A = -4 \left(-\frac{s}{2}+m +\frac{1}{2}\right) \left(s -1\right) {\mathrm e}^{-\frac{T \left(-1+m -s \right)}{2}}+m^{2} {\mathrm e}^{T}+m^{2} {\mathrm e}^{-T \left(m -s \right)}\\[6pt]
B = -8 \left(-\frac{s}{2}+m +\frac{1}{2}\right) \left(s -1\right) {\mathrm e}^{-\frac{T \left(-1+3 m -3 s \right)}{2}}\\[6pt]
C = \left(-2 m^{2}+\left(16 s -16\right) m -8 \left(s -1\right)^{2}\right) {\mathrm e}^{-T \left(-1+m -s \right)}\\[6pt]
D = -8 \left(-\frac{s}{2}+m +\frac{1}{2}\right) \left(s -1\right) {\mathrm e}^{-\frac{T \left(m -s -3\right)}{2}}\\[6pt]
E = \left({\mathrm e}^{-2 T \left(m -s \right)}+{\mathrm e}^{2 T}\right) m^{2}
\end{array}
$$
\begin{figure}
\center
\includegraphics[width=1\textwidth]{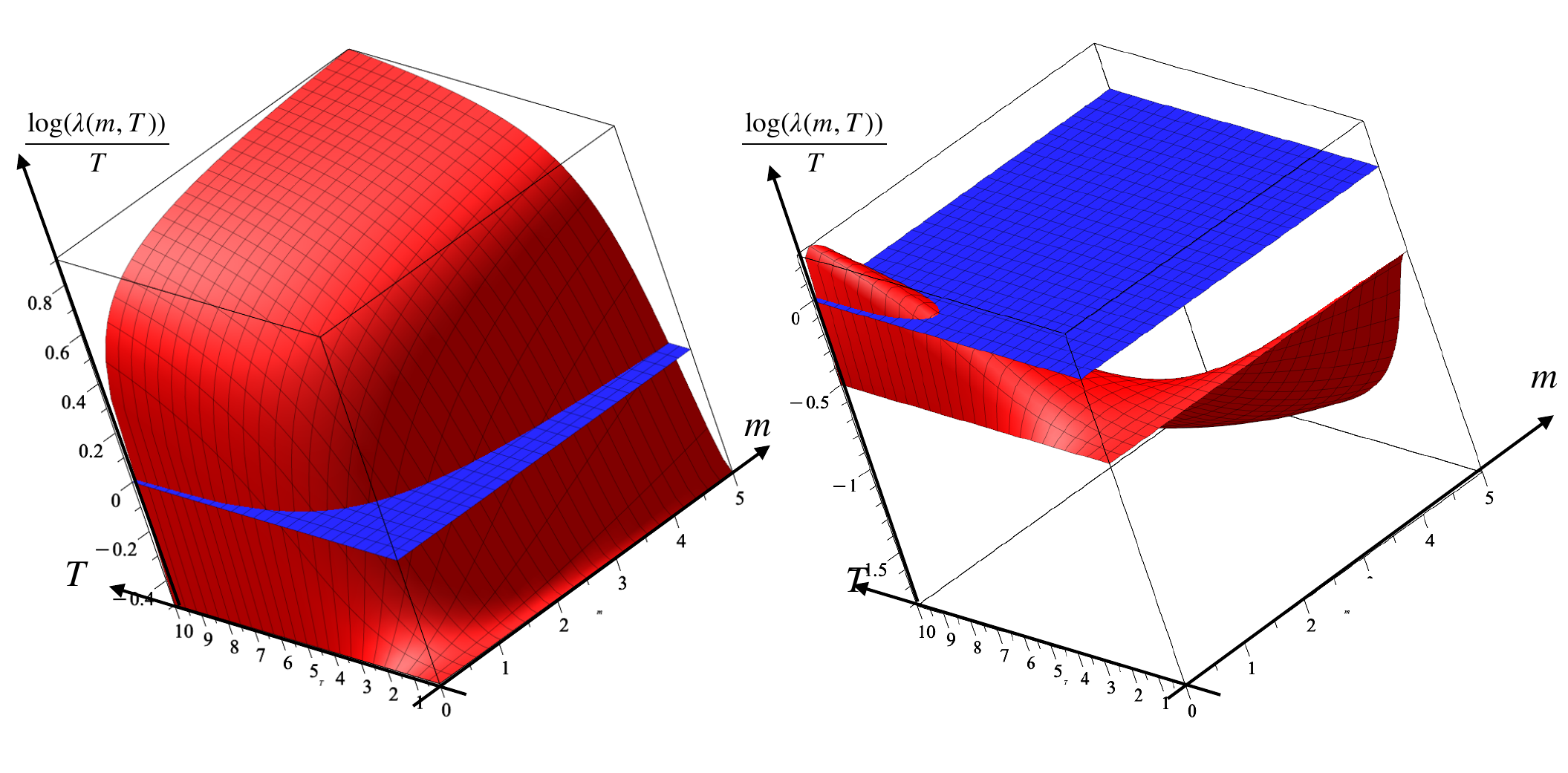}\\
\caption {On the left: growth ‘‘migration to the source'' \eqref{vlS}. 
On  the right   growth of system ‘‘migration to the sink''. Parameters:$ r_1^1=r_2^2 = -2;\,r_2^1= r_1^2 = 1$ }\label{Lyapcigognes}
\end{figure}
This is a complicated expression but we can see that when $m$ is large one can neglect all the terms which have $\e^{- Tm}$ in factor and it remains only
$$\lambda(m,T) \approx \frac{m^2 \e^{T} }{(1+m-s)^2}$$
 which is easily understood. 
 But we can also ask to Maple to  plot the exact graph of $(m,T) \mapsto \frac{\log(\lambda(m,T))}{T} $ which is done on Figure \ref{Lyapcigognes}-left. So the faster is the migration from the ‘‘sink'' to the ‘‘source'' te better will be the growth.
 
 This is entirely consistent with intuition since, ultimately, if populations are able to move instantly from one site to another, they will always benefit from the most favorable growth rate.
 
 %=======================
 \subsubsection*{Migration to the sink}

  It is interesting to look now to the same system but with the direction of migrations reversed as it is shown on Figure \ref{cig}-right i.e. the case of a population of storks having the curious idea to migrate to the worst site. 

 For this system  $\mbox{‘‘south'' } \to \mbox{ ‘‘north'' } \to \mbox{ ‘‘south'' } =\mathcal{C}_{sns}$  is a circuit such that $\chi^{\mathcal{C}_{sns}} > 0$ (actually it has the same growth index than previously). Thus, Proposition \ref{suffcroissance} applies, which is somewhat counter intuitive. How is it possible to have growth when the population systematically flees out the ‘‘source'' site to go to the ‘‘sink'' site? Like previously we can compute the monodromy matrix of this system and ask to Maple the expression of its dominant eigenvalue $\lambda(m,T)$ (that we do not show here) and the graph of $(m,T) \mapsto \frac{\log(\lambda (m,T))}{T}$ which is shown on the right of Figure \ref{Lyapcigognes}. This figure compares the Lyapunov exponents for the ‘‘migration towards the source'' and the ‘‘migration towards the sink'' for the same parameter value $r^i_j$. We can see that the two plots are very different. In the case of ‘‘migration towards the source'', except for the small blue region corresponding to low migration rates or to a short period, the growth rates are positive and increasing with $m$ and $T$.
On the contrary, in the case of ‘‘migration towards the sink'', the growth rate is essentially negative, which is expected, except for a small region corresponding to very  large values of $T$ and very  small values of $m$, which is not intuitive but was predicted by Proposition  \ref{suffcroissance}.

%====================================
 \subsection{More complex examples} \label{examples} 
  
  In the previous example, we saw that a positive index circuit, while predicting the possibility of DIG, does not characterize the growth mode at all. In this paragraph, we will go a little further by giving examples where, although all circuits have a negative growth index, there is still growth.

 %-----------------------------------------------
 \subsubsection*{ 3 sites, 2 seasons.}\label{32}
 \begin{minipage}{0.5\textwidth}
Consider the system $\Sigma $ defined by the scheme on the right. It is the $T$-periodic system $\Sigma(r,s,m,T)$ defined by
$$\frac{d \dr{x}}{dt}= A_1\dr{x} \quad \mathrm{if}\quad t\in \left[0,\frac{T}{2}\right)$$
$$\frac{d \dr{x}}{dt} =A_2\dr{x} \quad \mathrm{if}\quad t\in \left[\frac{T}{2},T\right)$$
with
\end{minipage}\quad
 \begin{minipage}{0.45\textwidth}
 \includegraphics[width=1\textwidth]{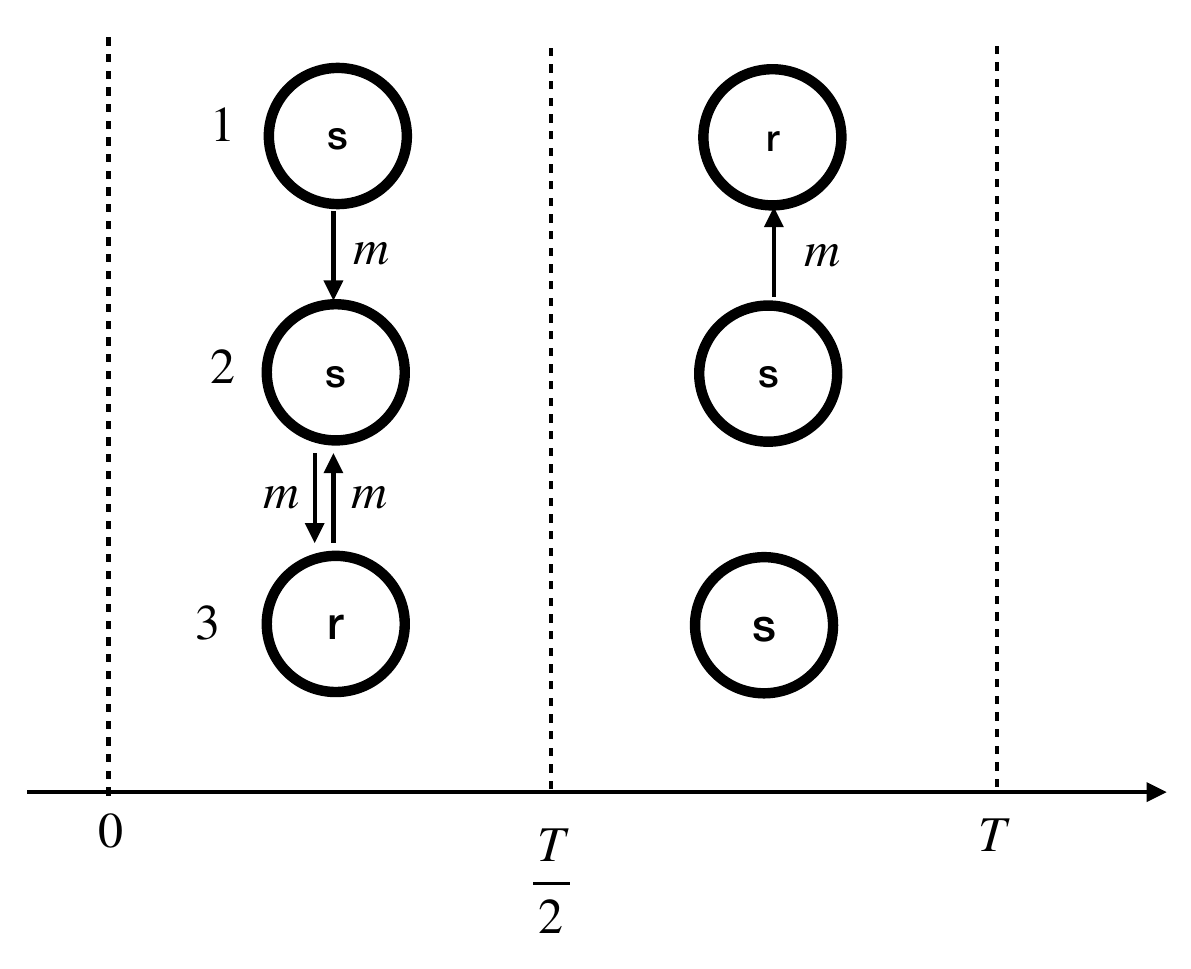}
\end{minipage}\\
\beq
A_1 = 
\left(
\begin{array}{ccc}
s-m&0&0\\
m&s-m&m\\
0&m&r-m
\end{array}
\right)
\quad \quad
A_2 =\left(
\begin{array}{ccc}
r&m&0\\
0&s-m&0\\
0&0&s
\end{array}
\right)\feq
We assume that $s <- r<0$. With this condition it is easily seen that, in the absence of migration, each site is a ‘‘sink''. As one can sea easily, the only 2-circuit of this system is the circuit $\mathcal{C} = |1\to 2||2\to 1|$ whose growth index $\chi^{\mathcal{C}}$ is $\frac{s}{2}+\frac{r}{2} <0$. Hence,
Proposition \ref{suffcroissance} does not apply and we cannot conclude to the existence of DIG. 

However, there are $T$ and $m$ such that the system is growing, as we will now show. Define 
\beq
M(r,s,m,T) = \e^{TA_1} \e^{TA_2}
\feq 
As we know the solution $\dr{x}(T)$  of $\Sigma(r,s,m,T)$ at time $T$, from the  initial condition $\dr{x}_0$ is given by
\beq
\dr{x}(T) = M(r,s,m,T)\dr{x}_0
\feq
Now we ask to Maple to compute the matrix $ M(r,s,m,T)$. With our computer Maple is not able to compute the eigenvalues of $ M(r,s,m,T)$ but fortunately it is sufficient to consider its entry $[M(r,s,m,T)]_{1,1}$ (first line, first column). Indeed, due to Proposition~\ref{orthaninvar} in Appendix~\ref{systlin}, we have 
$x_1(T) \geq [M(r,s,m,T)]_{1,1}x_1(0)$
and therefore
$x_1(mT) \geq [M(r,s,m,T)]_{1,1}^mx_1(0)$. Hence if $[M(r,s,m,T)]_{1,1} > 1$ the system is growing. From Maple we obtain

   $$
[M(r,s,m,T)]_{1,1} =  \frac{A+B+C+D+2E}{2\Delta\, \left(m -s +1\right) \left(s -1+\Delta \right) }
 $$
 with $\Delta = \sqrt{4m^2+(s-1)^2}$ and
  
 $$
 A = \left((1-s)^2 +\left(1-s\right) \Delta+2m^{2} \right) {\mathrm e}^{-\frac{T \left(-3+2 m -s +\Delta \right)}{4}}
 $$
$$
B =\left((1-s)^2-(1-s)\Delta+2m^2\right){\mathrm e}^{\frac{T \left(-4 m +3 s +1+\Delta \right)}{4}}
$$
$$
C = \left(-(1-s)^2-(1-s)\Delta-2m^2\right){\mathrm e}^{-\frac{T \left(4 m -3 s -1+\Delta \right)}{4}}
$$
 $$
 D = \left(-(1-s)^2 + (1-s)\Delta - 2m^2\right){\mathrm e}^{\frac{T \left(3-2 m +s +\Delta \right)}{4}}
 $$
 $$
 E= -\left(\left(m -2 s +2\right) {\mathrm e}^{-\frac{T}{2} \left(m -s -1\right)}+{\mathrm e}^{-T \left(m -s \right)} \left(s -1\right)\right) \Delta 
 $$
  \begin{minipage}{0.4\textwidth}
   We do not try to simplify nor analyse directly this expression but ask to Maple to plot the graph of $(m,T)\mapsto [M(r,s,m,T)]_{1,1}$. On the right, one sees the result in the case  $r = 1$ and $ s = -2$. As we observe, for sufficiently large values of $m$ and $T$ (for instance (5,10)) the entry $M(r,s,m,T)]_{1,1}$ is strictly greater than 1. Thus there is DIG which is not predicted by Proposition \ref{suffcroissance}.
\end{minipage}\quad
 \begin{minipage}{0.55\textwidth}
 \includegraphics[width=1\textwidth]{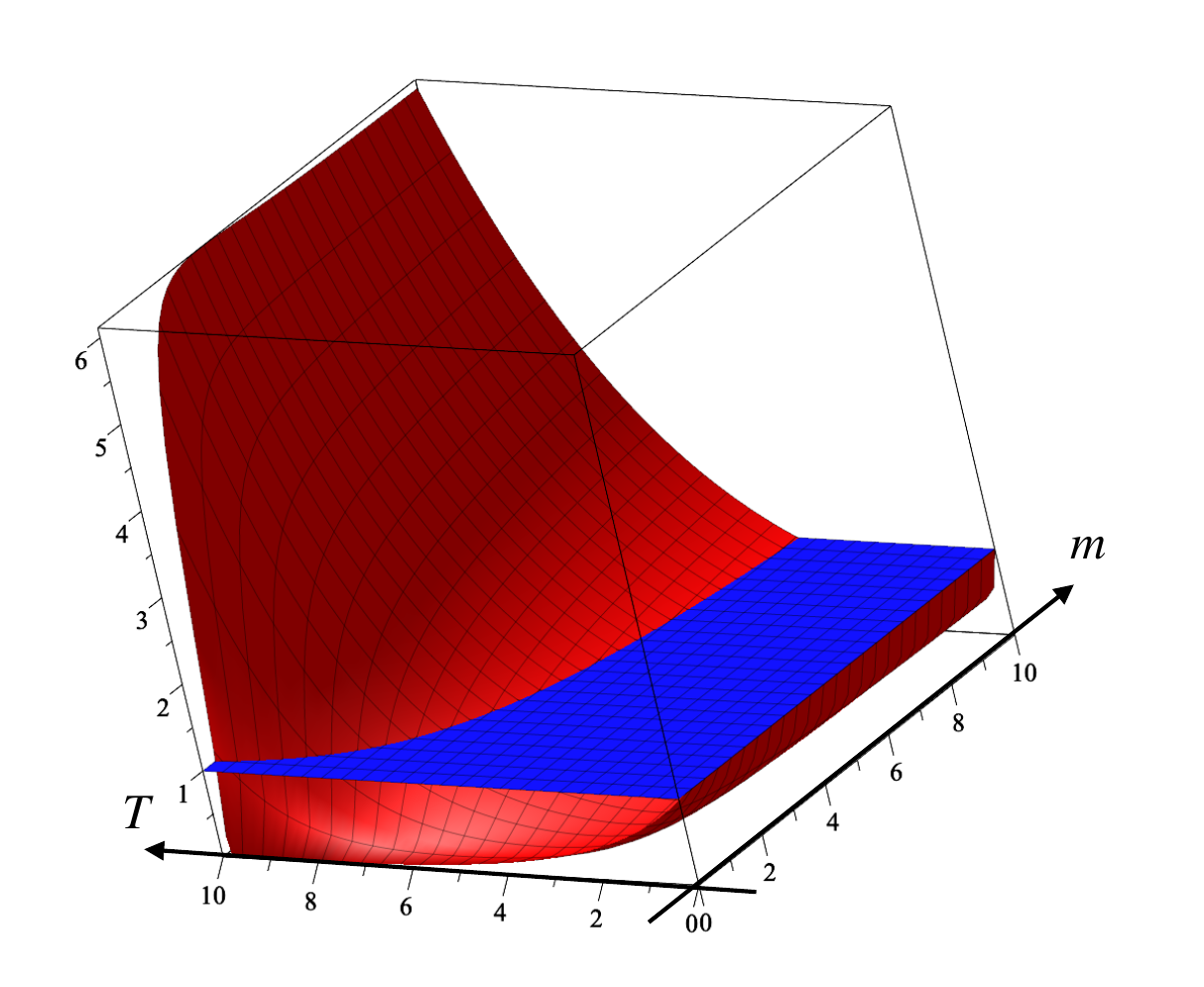}
\end{minipage}\\
If one looks carefully to this example one can see that the ‘‘circuit'' (which is not a circuit since it has a loop) $|1\to 2 \to 3 \to 2 ||2 \to 1|$ has an index which is $\frac{r}{2}+\frac{r}{2} = r$ which is strictly positive.

%---------------------------------------------------------
\subsubsection*{3 sites, 3 seasons}\label{33}

 \begin{minipage}{0.45\textwidth}
Consider the system $\Sigma $ defined by the scheme on the right. It is the $T$-periodic system $\Sigma(r,s,m,T)$ defined by\\
$  \displaystyle\frac{d \dr{x}}{dt}= A_1\dr{x} \quad \mathrm{if}\quad t\in \left[0,\frac{T}{3}\right)$\\
$ \displaystyle \frac{d \dr{x}}{dt} =A_2\dr{x} \quad \mathrm{if}\quad t\in \left[\frac{T}{3},\frac{2 T}{3}\right)$\\
$ \displaystyle \frac{d \dr{x}}{dt} =A_3\dr{x} \quad \mathrm{if}\quad t\in \left [\frac{2T}{3},T\right)$\\
\end{minipage}\quad
 \begin{minipage}{0.55\textwidth}
 \includegraphics[width=1\textwidth]{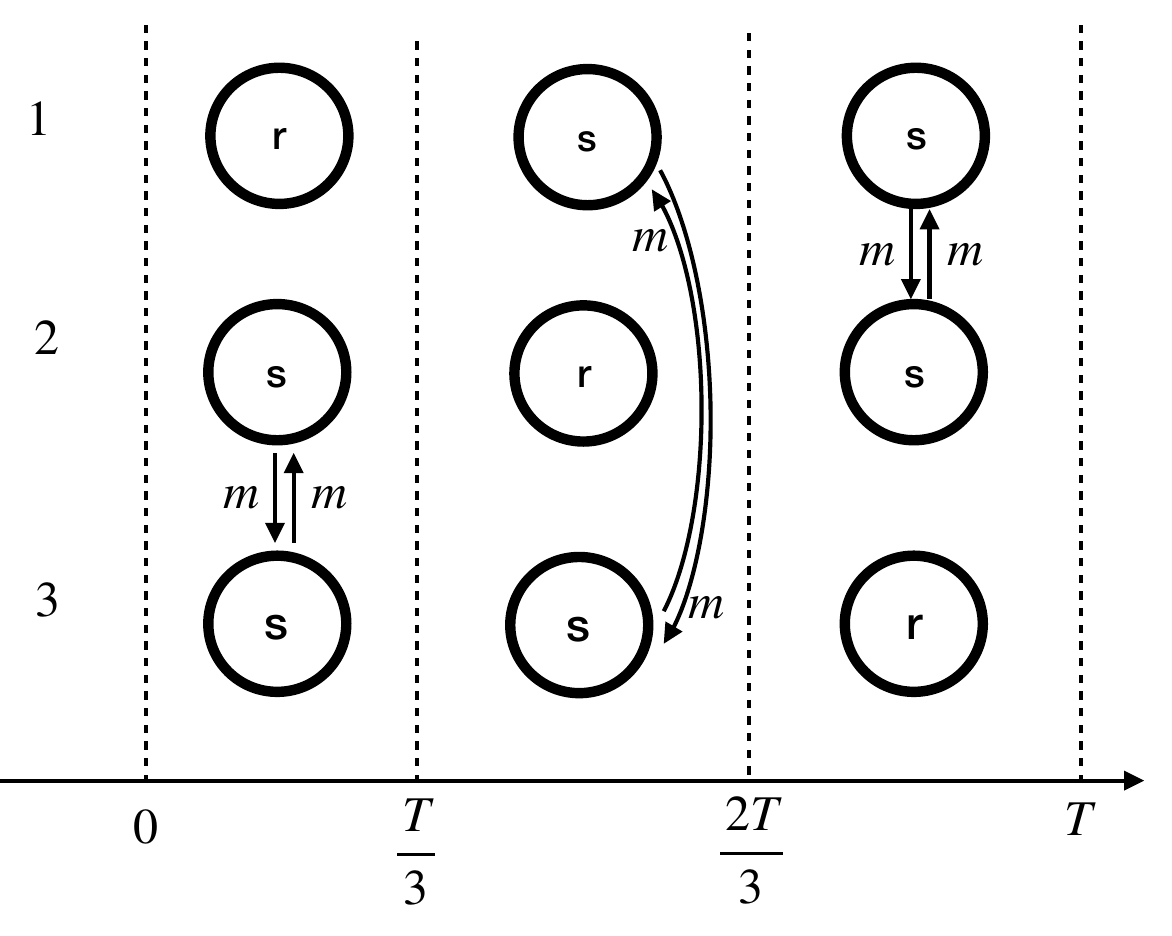}
\end{minipage}
 with
 \beq \label{M33}
\underbrace{
\left(
\begin{array}{ccc}r&0&0\\
0&s-m&m\\
0&m&s-m
\end{array}
\right)
}_{A_1}
\quad 
\underbrace{
\left(
\begin{array}{ccc}
s-m&0&m\\
0&r&0\\
m&0&s-m
\end{array}
\right)
}_{A_2}
\quad
\underbrace{
\left(
\begin{array}{ccc}
s-m&m&0\\
m&s-m&0\\
0&0&r
\end{array}
\right)
}_{A_3}
\feq
It was proved in \cite{BLSSJOMBa} (see § 4.5.2), thanks to the computation by Maple of the principal eigenvalue of $M = \e^{\frac{T}{3}A_3}\e^{\frac{T}{3}A_2}\e^{\frac{T}{3}A_1}$ that for $r = 1$ and $s = -1$ DIG is not present, while it occurs for $r = 1$ and $s = -0.8$ (see fig. 9 of \cite{BLSSJOMBa}).

Since, as it is readily seen, there is no $p$-circuit in this system, our Proposition \ref{suffcroissance} does not apply and thus is unable to predict the presence of DIG for 
$s = -0.8$.

The next Section is devoted to some extensions to Proposition \ref{suffcroissance} that enable us to deal with these two examples.
 %============================================
 \subsection{Better  sufficient condition for DIG }\label{BSCD}
  
 For simplicity of exposition we have stated and proved a simple growth condition (our Proposition \ref{suffcroissance}). But from this we can demonstrate the following more efficient condition which is suggested by the two previous examples.
 
 Let us consider the system \eqref{systemper2} with $p$ seasons.
 \begin{definition}$\;$
 Let us consider  a dynamic network $\mathcal{N}^{[1\cdots p]}$.
 \bito
 \item  A  \textbf{lpath} (path with loop):
 $$ a_0\to a_1 \to \cdots \to a_{i-1}\to a_i \to \cdots a_{l-1}\to a_l$$
 is a sequence of links connecting sites $a_1,\cdots,a_p$ \textbf{without requiring}, as in the case of  paths, that all the sites be different. 
 \item A \textbf{ dynamic $p$-lpath} is like previously  a sequence of $p$-lpaths which respect the seasons.
\item \textbf{ $q*p$-lcircuit}. Let $q$ be an integer. A  $q*p$-lcircuit $\mathcal{C}$ is a dynamic  $p$-lpath such that 
$$a_0 ≠ a_p ,\; a_0 ≠ a_{2p},...a_0 ≠ a_{(q-1)p},\; a_0 = a_{qp}$$
\item The \textbf{growth index} of a $q*p$-lcircuit $\mathcal{C}$ is:
\beq 
\chi{\mathcal{C}} = \sum_{k = 1}^{q\,p} (t_k - t_{k-1}) \max_{\pi_i \in \Gamma^k} r_i^k
\feq

\fit
\end{definition}

\begin{prop}\label{minocircuit2}
Consider the system $\Sigma(r_i^k ,\ell_{ij}^k,m,T)$ on the underlying network $\mathcal{N}^{[1\cdots p]}$. 
Consider a $q*p$-lcircuit $\mathcal{C}$ 
and its growth index $\chi^{\mathcal{C}}$. 
 Then there exist $\tau$ and constants $C >0$ and $\mu>0$ (independent of $m$) such that 
\beq \label{minoration2}
T > \tau \implic x_{i^1(0)}(T) \geq   C m^{\Lmat} \e^{T\left(\chi^{\mathcal{C}}- \mu \times m \right)} x_{i^1(0)}(0)
\feq
where $\Lmat = \Lmat(\mathcal{C})$ is  the total  length of the $q*p$-lcircuit. 
\end{prop}
\textbf{Proof.} The extension to a $q*p$-cirquit is straightforward since  it suffices to consider the system over a period $T' = qT$.
 In the Appendix \ref{extension} we explain how we can extend the lemma to a lpath with one loop. From the idea of this extension, it's not difficult to write a proof of the proposition. $\Box$

Similar to Proposition~\ref{suffcroissance} which is a corollary of Proposition~\ref{minocircuit}, we obtain from Proposition~\ref{minocircuit2} a better sufficient condition for DIG:

\begin{prop}
\textbf{{\em - Sufficient condition for DIG, 2}}.
 Consider the T-periodic system $\Sigma(r_i^k ,\ell_{ij}^k,m,T)$ defined by \eqref{systemper2}.  Assume that there exists a  $q*p$-lcircuit  $\mathcal{C}$ of the underlying time-varying network of $\Sigma(r_i^k ,\ell_{ij}^k,m,T)$ such that $\chi^\mathcal{C} >0$. Then for all $\tau$ large enough, there exist   $0 < a(\tau) < b(\tau)$ such that that for every $T \geq \tau $,   
$$ m\in[a(\tau),b(\tau)] \implic \;\Sigma(r_i^k ,\ell_{ij}^k,m,T)\mathrm{\; is \; a\;source}.$$
In other words, the existence of a $q*p$-lcircuit  $\mathcal{C}$  with $\chi^{\mathcal{C}} > 0$ is a sufficient condition for DIG.
\end{prop}
 Let's return to the example 3 sites , 2 seasons  of the Section  \ref{32}  and consider the 2-lcircuit  with a loop during season 1:
 	$$\underbrace{|1\to 2 \to 3 \to 2|}_{\mathrm{season 1}}\underbrace{|2\to 1|}_{\mathrm{season 2}}$$

\begin{figure}
\center
\includegraphics[width=1\textwidth]{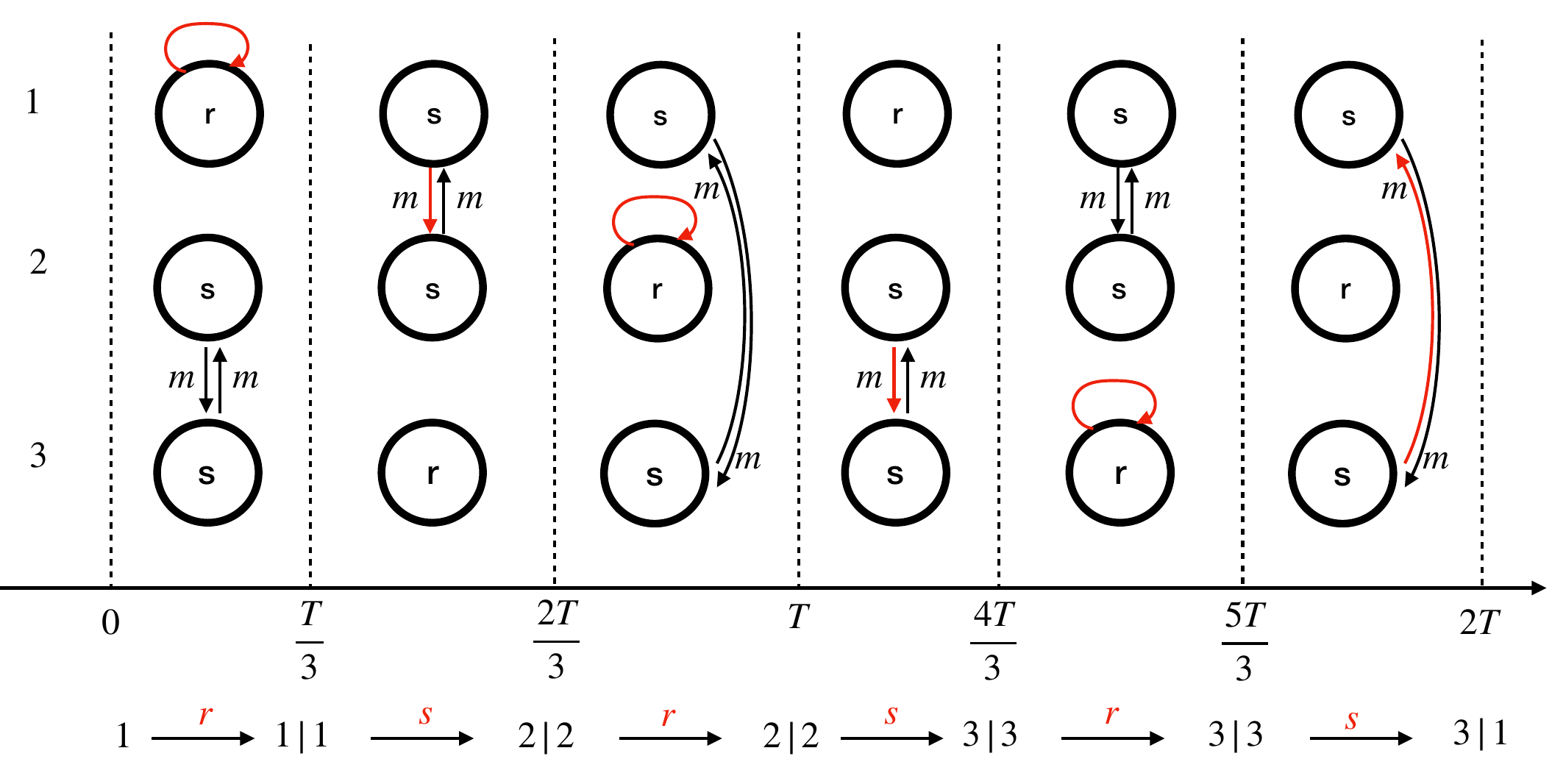}\\
\caption {Underlying network of system "3 sites,3 seasons" of Paragraph \ref{examples}  on a duration of 2-periods. Below the description of a 2*3-lcircuit.}\label{3S3S}
\end{figure}

\noindent The growth index of such a circuit, calculated as before by taking the dominant growth rate on consecutive paths, is:
	$$\frac{1}{2} r + \frac{1}{2}r = r$$
which is positive and Proposition \ref{minocircuit2} predicts the growth.\\\\
 Let's now take the example of 3 sites over 3 seasons given  in Paragraph \ref{examples} and 
observe a duration of 2 periods as shown in Figure \ref{3S3S}.
On this figure we can observe the dynamic $2*3$-lcicuit
\beq
|1\to 1||1 \to 2||2\to 2||2 \to 3||3 \to 3||3 \to ||3 \to 1|
\feq 
and remark that the succession of dominant growth rates on the six successive paths is
\beq
r||s||r||s||r||s
\feq
and the growth index is $\frac{1}{3} r+\frac{1}{3} s+\frac{1}{3} r+\frac{1}{3} s+\frac{1}{3} r+\frac{1}{3} s = r+s$ which is positive as soon as $r > -s$. Thus Proposition \ref{minocircuit2} predicts, for instance, that 
$$ r = 1 ,\; s > - 1 \implic \mathrm{growth\;for\;suitable}\;(m,T)$$

\begin{figure}
\center
\includegraphics[width=0.7\textwidth]{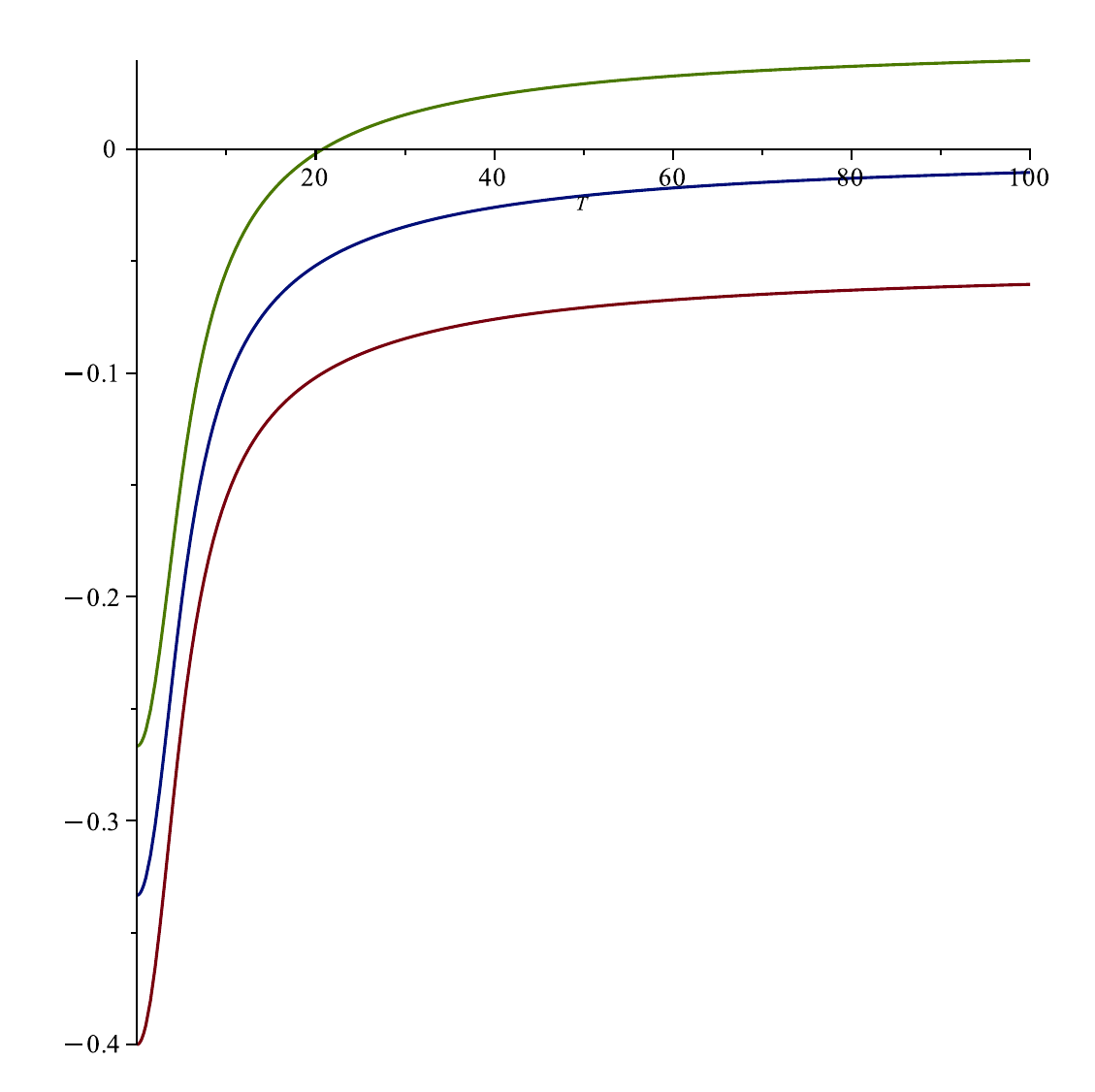}\\
\caption {Graph of $T \mapsto   \Lambda(0.5, T)$ (see \eqref{lyap}) for $s= - 0.9$ (green), $s = -1$ (blue), $s = -1.1$}\label{Ly3S3S}
\end{figure}

The Lyapunov exponent of system "3 sites,3 seasons":
\beq \label{lyap}
\Lambda(m,T) = \frac{\log(\lambda(m,T))}{T}
\feq
where $\lambda(m,T)$ is the dominant eigenvalue of $\e^{\frac{T}{3}A_3}\e^{\frac{T}{3}A_2}\e^{\frac{T}{3}A1)}$, gives the asymptotic  growth rate of the system. It can be computed by Maple and we display on Figure \ref{Ly3S3S} the graph of $T \mapsto \Lambda(0.5, T)$. We see that with $s = -0.9$  there is growth as soon as $T> 20$.

%====================================================
\section{Extension to random perturbations of the duration of seasons}\label{stochastique}
%====================================================

In this section, we study the case where the duration of each season is random, instead of being deterministic. We prove that we can define a mean growth index associated to a $p$ - circuit, and that if one can find a $p$-circuit with a positive mean growth index, then there is DIG and the growth threshold is exponentially small.

 More precisely we consider a succession of ‘‘years'' (or ‘‘cycles'') indexed by $j$, composed of a succession of $p$ ‘‘seasons'' indexed by $k$. The $k$-th season of the $j$-th year has a random duration $TU^{j,k}$, where $T >0$ is a parameter and $\mathbf{U}^j = (U^{j,k})_{1 \leq k \leq p} $ is a random vector  with values in $\mathbb{R}_+^p$. We let $t_0^j =t_k^0= 0$ and $t_k^j = t_p^{j-1}+\sum_{l=1}^{k} U^{j,l}$. Then, for $j \geq 1$, $t_k^j$ denote the time in year $j$ where the season is switched from $k$ to $k+1$, while $t_p^j-t_p^{j-1}$ denote the length of year $j$.  This defines the system\\\\
 \fbox{
\begin{minipage}{0.98\textwidth}
\beq \label{systemrandom}
\Sigma(r_i^k ,\ell_{ij}^k,\mathbf{U}^j,m,T)=\quad \quad\left\{
\begin{array}{l}
\displaystyle \frac{dx_i}{dt} = r_i^kx_i+ m \sum_{j=1}^n \ell_{ij}^k x_j  \quad t \in [Tt^j_{k-1},\,Tt^j_k)\\[8pt]
i = 1,...,n\quad\quad j \geq 1, \quad\quad k = 1,...,p 
\end{array}
\right.
\feq
$$\;$$
\end{minipage}}\\\\

%\begin{exemple}
%A natural example is when we assume that there is a small perturbation in the \textit{date} at which the season is switched. In the periodic case studied in the previous section, the system goes to the $k$-th season to the $(k+1)$-th at time $t_k$. Let $\delta \leq \min_k | t_k - t_{k-1}|$ be the minimal  deterministic length of a season, and let $(u_k^j)_{k=1,\ldots,p}$ be a family of independent variables such that $|t_k - u_k^j| < \delta$ for all $k=1, \ldots, p$, and we assume that during year $j$, the passage from season $k$ to $(k+1)$ occurs at time $t_k + u_k^j$. If $\delta$ is small, this can be seen as a small random perturbation of the arrival time of the $k$-th season. This is encompassed in our framework by setting 
%\end{exemple}

We make the following assumption on the length of the season, and the independence of the length from year to year:

\begin{hyp}
\label{Hrandom}
For all $k = 1, \ldots, p$, $\mathbb{P}( U^{j,k} = 0 ) = 0$ and $\mathbb{E}( U^{j,k}) =\mathbb{E}( U^{1,k})   < + \infty$. In addition, we assume that $(\mathbf{U}^j)_{j \geq 1}$ is a sequence of independent and identically distributed random variables.
\end{hyp}
\begin{rem}
Note that we do not require the independence of the $(U^{j,k})_{k=1, \ldots, p}$ for fixed $j$. In particular, one may for example assume that $t_p^j - t_p^{j-1}=\sum_{k=1}^p U^{j,k} =1$, so that every year has the same length. This holds in the periodic case, corresponding to $U^{j,k}  = t_k - t_{k-1}$ almost surely, for all $j \geq 1$ and $k = 1, \ldots, p$.
\end{rem}
%\red{Dans le cas périodique la somme des  $t_k - t_{k-1}$ vaut $1$ ce que tu ne semble pas supposer pour tes v.a. Tu ne dis rien sur leur loi. A mon avis ce qui fait un peu sens en termes de modélisation c'est 1) le choix d'une durée totale $T_{aleat}$ qui pourrait être quelque chose de centré autour d'une moyenne $T$ qui serait le paramètre ; 2) un découpage aléatoire (mais pa trop) autour de valeurs moyennes. Pour moi un truc simple qui ferait sens serait de prendre des intervalles définis par des bornes $t_j + u_j$ où les $u_i$ sont aléatoires et bornées de façon à ce que les $t_j$ restent bien ordonnées et ensuite on multiplie tout ça par $T$. Une fois de plus on retrouve le fait que prendre la période comme paramètre dans le modèle déterministe n'est pas une bonne idée.}

The definition of a $p$-circuit is the same as in the periodic case and  we precise the definition of the growth index in that case:

\begin{definition} \textbf{Growth index of a circuit with random duration} 

Let 
\beq \label{pcircuit-random}
\mathcal{C} =
a_0\Gamma^1 a_1\Gamma^2 a_2 \cdots a_{k-1}\Gamma^k a^k \cdots
a_{p-1}\Gamma^pa_p
\feq
with 
$$ a_{k-1}\Gamma^ka_k = \left\{  a_{k-1} = \pi_{i^k(0)} \to \pi_{i^k(1)}\to \cdots \pi_{i^k(j)} \to \pi_{i^k(l^k)}= a_k\right \}$$
 be a  $p$-circuit defined on the underlying dynamic network of the system \eqref{systemper3}.  We call {\em random growth index in the} $j${\em -th year} of the p-circuit $\mathcal{C}$ the number
\beq 
\chi_{stoc}^{j,\mathcal{C}} = \sum_{k = 1}^p U^{j,k} \max_{j = 0}^{l^k} r^k_{i^k(j)}
\feq
and {\em mean growth index} (of the $p$-circuit $\mathcal{C}$) is the number 
\beq 
\chi_{stoc} ^{\mathcal{C}} = \sum_{k = 1}^p \mathbb{E}(U^{1,k}) \max_{j = 0}^{l^k} r^k_{i^k(j)}
\feq
\end{definition}

\begin{rem}
Note that in the above definition, the $p$-circuit $\mathcal{C}$ is deterministic, and is the same year after year. The random growth index $\chi_{stoc}^{j,\mathcal{C}}$ in the $j$-th year of $\mathcal{C}$ corresponds to the value that would have the growth index of $\mathcal{C}$ in the periodic case if $t_{k} - t_{k-1}$ were set equal to $U^{j,k}$ every year. The mean growth index $\chi_{stoc} ^{\mathcal{C}}$ as for it, is the value that would have the growth index of $\mathcal{C}$ if $t_{k} - t_{k-1}$ were set equal to the mean duration of the season $k$, that is $\mathbb{E}(U^{j,k})$.
\end{rem}

Using exactly the same proof as for Proposition~\ref{minocircuit}, we  can prove:
\begin{prop}\label{minocircuitrandom}
Consider the system \eqref{systemper3}  on the undelying network \eqref{UN}. Consider the $p$-circuit $\mathcal{C}$ defined by \eqref{pcircuit-random} and its random growth indexes $(\chi_{j,stoc}^{\mathcal{C}})_{j \geq 1}$. Then for all $T > 0$ and all $j \geq 1$,
\beq
  x_{i^1(0)}(T(U^{j,1} + \ldots + U^{j,p})) \geq   \left( \prod_{k=1}^p C_k(T U^{j,k}) \right) m^{\Lmat} \e^{T(\chi_{stoc}^{j,\mathcal{C}}- \mu m)} x_{i^1(0)}(T(U^{j-1,1} + \ldots + U^{j-1,p}))
\feq
where ${\Lmat }$ is the total length of the circuit, $C_k(\cdot)$ is the increasing function given by Lemma \ref{lemme} (in Appendix~\ref{integrationchemin}) on the path $\pi_{i^k(0)}\Gamma^k \pi_{i^k(l_k)}$, and $U^{j-1,k}=0$ for all $k=1, \ldots, p$.
\end{prop}

\begin{rem}
The proof of this Proposition is obtained as follows: for a given year $j$ and a realization of $U^{j,1}, \ldots, U^{j,p}$, we apply Proposition~\ref{minocircuit} to the deterministic periodic system with $t_k - t_{k-1} = U^{j,k}$. Note that the bound is random, since although the functions $C_k$ are deterministic, they are evaluated at a random value $TU^{j,k}$.
\end{rem}

Now, as in the deterministic case, we  introduce the DIG threshold as 
\[
m^*(T) = \inf_{m > 0} \{ m: \;  S(t,m,T) \to \infty \},
\]
where $S(t,m,T)$ is the total population of $\Sigma(r_i^k ,\ell_{ij}^k,\mathbf{U}^j,m,T)$.
\begin{rem}
Without additional assumptions, $m^*(T)$ is a priori random. However, in the case where the total migration matrix is irreducible (see Remark \ref{rem:totalmigration-Lambda}), it can be proven as in \cite{BLSSJOMBb} and \cite{BLSSAFST}  that there is a \textbf{deterministic} Lyaponuov exponent $\Lambda(m,T)$ such that, for any initial condition, almost surely, 
\[
\Lambda(m,T) = \lim_{t \to \infty} \frac{1}{t} \log( S(t,m,T) ) .
\]
This implies that $m^*(T)$ is almost surely constant in that case.
\end{rem}

Using Proposition~\ref{minocircuitrandom}  and the law of large numbers, we  now prove that the DIG threshold is exponentially small with respect to the parameter $T$:

\begin{prop}
\label{prop:exponentiallysmallrandom}
Assume that $\mathbb{E}(|\log(U^{j,k})|) < + \infty$ and that there exists a circuit $\mathcal{C}$ with $\chi_{stoc}^{\mathcal{C}} > 0$. Then, for all $1 > \varepsilon >0$, there exists $\tau > 0$ such that for all $T \geq \tau$, one has almost surely,
\[
m^*(T) \leq e^{ -\frac{1}{{\Lmat }}(1 - \varepsilon)T \chi_{stoc}^{\mathcal{C}}}
\]
\end{prop}

\noindent \textbf{Proof.}
Without loss of generality, we assume that $i_1(0) = 1$. Set $X^0 = x_1(0)$ and for all $n \geq 1$, 
\[
X^n = x_1\left( T \sum_{j=1}^n \sum_{k=1}^p U^{j,k}\right),
\]
the size of the population on site $1$ after $n$ years. 
Then, Proposition~\ref{minocircuitrandom} implies that
\[
X^n \geq \prod_{j=1}^n Y^j \, X^0,
\]
where  for all $j=1, \ldots, n$,
\[
Y^j = \left( \prod_{k=1}^p C_k(T U^{j,k}) \right) m^{\Lmat } \e^{T(\chi_{stoc}^{j,\mathcal{C}}- \mu m)}
\]
Note that $Y^j$ is a sequence of i.i.d. random variables, such that $\mathbb{E}(|\log Y^j|) < + \infty$. Hence, the strong law of large numbers implies that, almost surely, 
\[
\lim_{n \to \infty} \frac{1}{n} \sum_{j=1}^n \log(Y^n) = \mathbb{E}( \log Y^1).
\]
Therefore,
\[
\limsup_{n\to \infty} \frac{1}{n} \log( X^n) \geq \mathbb{E}( \log Y^1),
\]
and the system is growing provided $\mathbb{E}( \log Y^1) > 0$. Let 
\[
\tilde C(T) = \sum_{k=1}^p \mathbb{E}\left(\log C_k(T U^{1,k}) \right),
\]
then 
\[
\mathbb{E}( \log Y^1)  = \tilde{C}(T) + {\Lmat } \log(m) + T (\chi_{stoc}^{\mathcal{C}} - \mu m)
\]
Let $\varepsilon > 0$. We can assume that $m < \frac{\varepsilon}{2} \frac{\chi_{stoc}^{\mathcal{C}} }{\mu }$, so that
\[
\mathbb{E}( \log Y^1)  \geq \tilde{C}(T) + {\Lmat } \log(m) + T (1 - \frac{\varepsilon}{2}) \chi_{stoc}^{\mathcal{C}}
\]
Now, it is easily seen from the definition of $C_k$ in Lemma~\ref{lemme} and the fact that $U^{1,k} > 0$, that for all $k = 1, \ldots, p$, $\log(C_k(T U^{1,k})$ converges monotonically to a finite limit as $T$ goes to infinity. Hence, by monotone convergence, $\frac{1}{T}\tilde{C}(T) \to 0$ as $T$ goes to infinity. Therefore, for some $\tau > 0$ and all $T \geq \tau$, we have $ \tilde{C}(T) \geq - \frac{\varepsilon}{2} T \chi_{stoc}^{\mathcal{C}}$, and thus
\[
\mathbb{E}( \log Y^1)  \geq  {\Lmat } \log(m) + T (1 - \varepsilon) \chi_{stoc}^{\mathcal{C}}.
\]
In particular, $\mathbb{E}(\log Y^1) > 0$ whenever $\frac{\varepsilon}{2} \frac{\chi_{stoc}^{\mathcal{C}} }{\mu } > m \geq e^{ -\frac{1}{{\Lmat}}(1 - \varepsilon)T \chi_{stoc}^{\mathcal{C}}}$, which concludes the proof of the proposition. $\Box$

\begin{rem}
Note that we have chosen to take the same $p$-circuit, independently of the realization of $U^{j,k}$. However, there might be situations where it would be more clever to adapt the $p$ - circuit from year to year, depending on the realization of the $U^{j,k}$. This is left for future work.
\end{rem}

%\red{Je ne sais pas si le reviewer 2 avait des faibles compétences en proba comme moi. Je ne sais pas ce que pensera F. Bienvenu qui lui a des compétences. Les reviewers avaient l'air satisfaits de la partie déterministe. On pourrait peut être la supprimer ou se contenter de donner des indications sur des futures investigations ?}

%===============
\section{Discussion}
%================
In this article, we considered the evolution of populations at different sites linked by migration paths. Environmental conditions vary periodically over time, as do migration rates between each site. The models are  continuous-time models.
We wanted to highlight how the evolution of the structure of the migration network influences total population growth. But rather than give precise growth results for a given model, we have sought to identify a method for minimizing growth, based on certain properties of the dynamic graph underlying the dynamic system.

To do this, we defined the growth index $\chi^{\mathcal{C}} $ of a simple $p$-circuit $\mathcal{C}$. A simple $p$-circuit is a route from site to site that respects the migration links existing during a given season and returning to the starting point. Our main result, Proposition \ref{minocircuit} is that {\em  if there is a simple $p$-circuit with a strictly positive growth index, then the total population is growing for some values of $m$ and $T$.} Following \cite{Katriel} we called this the DIG (Dispersal Induced Growth) effect.

Perhaps the most  important point that we can emphasize here is that  we do not make the assumption that migration matrices are irreducible,  which allows us to cover realistic cases such as seasonal population migrations from one site to another.

To keep things  mathematically as simple as possible - we're only using very elementary mathematical results from undergraduate courses  - we've chosen
to consider only periodic systems which coefficients are piecewise constant. There is little doubt that they are still true in the piecewise continuous case but this deserves further investigations.

One of the advantages of the ‘‘piecewise constant'' hypothesis is that it can be immediately extended to systems where  the duration of the seasons is no longer a fixed quantity but a random one, which is obviously much more realistic. As an example of what can be done, in Section \ref{stochastique} we define a stochastic version of the growth index along a circuit and demonstrate that the growth threshold is exponentially small with respect to the parameter $T$. 

In models of population dynamics forced by a periodic environment, the period is often fixed, e.g. year, month, and it is not very relevant, apart from mathematical considerations, to take, as we did,  the period as a parameter. A statement like “If the duration of the year is long enough, then...” doesn't make much sense. For instance, a model of the form
$$\Sigma(r_i(.),\ell_{ij}(.), m, S)\quad \quad \frac{dx_i(t)}{dt} =\big( r_i(t)+ S \sigma_i(t) \big)x_i(t) + m \sum_{j= 1}^n \ell_{ij}(t) x_j$$
with $r_i(\cdot), \sigma_i(\cdot), \ell_{ij}(\cdot)$ have a fixed period $1$,  where $S\sigma_i(\cdot)$ is a fluctuation around some average value $r_i$  is more relevant to express modification of environmental parameters as a function, for instance, of altitude, or latitude. The parameter $S$ which measures the amplitude of the fluctuations  replaces $T$ and seems more adequate for interpretation of mathematical results. There is no doubt that our method of minorization along paths also works  for $\Sigma(r_i(\cdot),\ell_{ij}(\cdot), m, S)$  but it remains to be done.

One area where our approach seems promising is that of epidemiology. Indeed, \cite{HOLTPNAS20,HOLT23} have shown that the phenomenon of inflation plays a negative role in the persistence of infected subjects. Network models are also well developed in this field and the analysis of the T-circuits as defined here could be useful by detecting where suppressing contact between infected and susceptible people is most effective.

There is a large body of literature (see \cite{SCH23} and its bibliography for a recent example) on the question of inflation in discrete-time models. 
Insofar as in a discrete-time model it is possible to transfer the entire population instantaneously from one site to another, which is not possible for continuous-time models, the questions that arise in the two cases are not exactly the same and, as a result, the comparison of results is not immediate. This could be the subject of further work.

Finally, we must make the following observation. On examples with few sites like the ones we have looked at, it is not difficult to determine the circuits and therefore calculate the associated growth indices. But what about systems with a large number of sites? This is a question of graph theory that we have not yet addressed.

%%%%%%%%%%%%%%%%%%%%%%%%%%%%%%%%%%%%%%%%%%%%%%%%%%%%%%%%%
\appendix 
%=====================================
 \section{Linear differential equations.}\label{systlin}
 %=====================================
 
 For the convenience of the reader we recall some elementary facts regarding linear differential equations and systems that can be found in elementary textbooks.
 
 %==============================================================
 \subsection{Closed form solutions for non autonomous linear differential equations }\label{eqlineaire}

Let $t \mapsto a(t)$ and $t \mapsto b(t)$ be two integrable functions. 
\begin{prop}\label{propfond}
  Let  $x(t,x_0)$  be the solution of the initial value problem 
\beq
\frac{dx}{dt} = a(t) x+b(t)\quad \quad x(0) = x_0
\feq
If one denotes  $ \displaystyle m(t) = \int_0^ta(\tau)d\tau$ one has
\beq \label{formulesol}
x(t,x_0) = \e^{m(t)}\left(x_0+\int_0^t \e^{-m(s)}b(s)ds\right).
\feq 
\end{prop}
If $a(t)$ is just a constant, \eqref{formulesol} reads 
\beq \label{formulesol2}
x(t,x_0) = \e^{t\,a}\left(x_0+\int_0^t \e^{-sa} b(s)ds\right).
\feq
 \subsection{Linear systems.} 
 \subsubsection*{Notations}
We use the following notation: for $\dr{x},\,\dr{y}  \in \Rmat^N$, $\dr{x} ≥ \dr{y} $ means that for all $i$, $x_i \geq y_i $ ; $\dr{x} > \dr{y} $ means that $x_i \geq y_i $ and $\dr{x} \not = \dr{y} $ ; and $\dr{x} \gg \dr{y}  $ means that for all $i$, $ x_ i >y_i $. We use the same notation for $n \times p$ matrices considered as elements of $\Rmat^{n\times p}$. Given the  system of differential equations in $\Rmat^n$  
$$ \Sigma \quad \quad \left\{ \frac{d\dr{x}}{dt} = f(\dr{x},t)\right. $$
we denote $\dr{x}(t,(\dr{x}_0,t_0))$ its solution  with initial condition $\dr{x}(t_0) = \dr{x}_0$. We say that $\Sigma$ is positive ($\Sigma \geq 0$) if $\dr{x}_0 \geq 0 \Longrightarrow \dr{x}(t,(\dr{x}_0,t_0)) \geq 0$ and given two positive systems 
$$ \Sigma_1 \quad \quad \left\{ \frac{d\dr{x}_1}{dt} = f(\dr{x}_1,t)\right.\quad \quad \Sigma_2 \quad \quad \left\{ \frac{d\dr{x}_2}{dt} = f(\dr{x}_2,t)\right. $$
we say that $\Sigma_1$ minorizes $\Sigma_2$ (denoted $\Sigma_1 \leq \Sigma_2$)  if
$$ 0 \leq \dr{x}_{1_0} \leq \dr{x}_{2_0}  \Longrightarrow  \dr{x}_1(t,(\dr{x}_{1_0},t_0)) \leq \dr{x}_2(t,(\dr{x}_{2_0},t_0))$$

%===============================
 \subsubsection*{Exponential of a matrix.}
  The exponential of a matrix allows one to extend formula \eqref{formulesol2} to linear systems.
 Let $\dr{x} = (x_1,\cdots,x_i,\cdots x_n)$ be a vector of $\Rmat^n$ and $A$  an $n \times n$ matrix.
 
  The matrix 
 $$ \mathrm{Id} + t  A + \frac{t^2}{2 !} A^2 +\cdots + \frac{t^k}{k !}A^k+ \cdots = \sum_{k = 0}^{\infty} \frac{t^k}{k!}A^k$$
 where $\mathrm{Id}$ is the identity matrix, is well defined (the sum is convergent) for every values (positive or negative) of $t$ and is denoted $\e^{tA}$. It has the following properties:
 
 \bito
 \item $t \mapsto \e^{t\,A}\dr{x}_0$ is the solution of the differential equation
 $$\frac{d \dr{x}}{dt} = A \dr{x}\quad \quad \dr{x}(0) = \dr{x}_0$$
 \item For every $t_1$and every $t_2$ one has $ \e^{t_1A} \e^{t_2A} = \e^{(t_1+t_2)A }$
 \item For every $t$ the matrix  $\e^{tA}$ is invertible and $\left(\e^{tA}\right)^{-1} = \e^{-tA}$ 
\item Il $t \mapsto \dr{b}(t)$ is an integrable mapping from $\Rmat$ to $\Rmat^n$ the solution of 

\beq \label{nonhomo}
\frac{d \dr{x}}{dt} = A \dr{x}+\dr{b}(t)\quad \quad \dr{x}(0) = \dr{x}_0
\feq
 is given by
 \beq \label{solnonhomo}
 \dr{x}(t) = \e^{tA} \dr{x}_0 + \int_0^t \e^{(t-\tau )A} \dr{b}(\tau) d\tau
 \feq
 \fit
 
 \begin{prop}\label{anulation}
 Let $\dr{x}_0 \not = 0$. Then, for every $t$ one has $\e^{tA}\dr{x}_0 \not = 0$ 
 \end{prop}
 \textbf{Proof}. This follows from the fact that  $\e^{tA}$ is invertible.$\Box$
%__________________________________________________________________
 \subsubsection*{Invariance of the positive orthant for Metzler systems.}\label{metzler}
 A matrix $M = \left(m_{ij}\right)$ is a Metzler matrix if $i\not = j \implic m_{ij} \geq 0$.  \begin{prop} \label{orthaninvar}
 If $M$ is a Metzler matrix then for every $\dr{x}_0 \geq 0$ and every $t \geq 0$ one has $\dr{x}(t) = \e^{tM}\dr{x}_0 \geq 0$. In particular, all the entries of $e^{tM}$ are nonnegative.
 \end{prop}
 \textbf{Proof}. If $\dr{x}_0 = 0$ one has $\e^{tM}\dr{x}_0 \equiv 0$ and the proposition is proved. We assume $\dr{x}_0 \not = 0$.\\
 {\em First step}. We assume that
 \beq \label{strict} 
 \forall \, i \;\; {x_i}_0> 0 \;\mathrm{ and }\; \forall \,i\;\forall_j: i\not = j \implic m_{ij} >0
 \feq
  and suppose that $e^{tM}\dr{x}_0$  is not positive for every $t > 0$. From Proposition \ref{anulation} it is not possible that all the components $x_i(t)$ vanish at the same time and thus, if all the components $x_i(t)$ are not always strictly positive,  there must exists $t^* > 0$ (the instant when a component vanishes for the firs time)  with the following properties:
 \ben
 \item There exists $i$ such that $x_i(t^*) = 0$ and $0\leq t < t^* \implic x_i(t) > 0$
 \item There exists at least one  $j \not = i$ such that $x_j > 0$
 \fen
 Since $\dr{x}(t)$ is a solution of $\displaystyle \frac{d \dr{x}}{dt} = M \dr{x}$ one has
 $\displaystyle  \frac{dx_i}{dt} = m_{ii} x_i(t) + \sum_{j \not = i} m_{ij}x_j(t)$ and, from 2) above, for $t = t^*$  
 $$  \frac{dx_i}{dt}(t^*) = m_{ii} x_i(t^*) + \sum_{j \not = i} m_{ij}x_j(t^*) = \sum_{j \not = i} m_{ij}x_j(t^*) > 0$$
 which contradicts point 1). As a consequence such a $t^*$ cannot exist and we have proved that under hypothesis \eqref{strict} 
 \beq \label{conclusionstrict}
 \forall \; t, \forall\;i\quad x_i(t) > 0
 \feq
 {\em Second step.} Let $M$ be a Metzler matrix and $\dr{x}_0 \geq 0$. Let
 $$ M_k = M + \frac{1}{k} \left(
 \begin{array}{ccccc}
 0&1&1&\cdots&1\\
 1&0&1&\cdots&1\\
 \cdot&\cdot&\cdot&\cdots&\cdot\\
 1&\cdots&1&0&1\\
 1&\cdots&1&1&0
 \end{array}
 \right)
 \quad \quad \dr{x}_k = \dr{x}_0 + \frac{1}{k}  \left(
  \begin{array}{c}
 1\\
 1\\
 \cdots\\
 1\\
 1
 \end{array}
 \right)
 $$
 We know that for every $t\geq 0$ one has 
 $$\lim_{k \to \infty} \e^{t M_k}\dr{x}_k = \e^{t M}\dr{x}_0 $$ 
 Since $M_k, \;\dr{x}_k$ satisfy \eqref{strict}, each component  $\dr{x_i}_k(t)$ is strictly positive and its limit is positive or equal to $0$, which proves the proposition. $\Box$.
 %----------------------------------------------------
 \subsubsection*{Comparison of solutions}
 We prove the following  proposition which is about the comparison of solutions of Metzler systems.
 \begin{prop}
 Let $M = \left(m_{ij}\right) $ and $N= \left(n_{ij}\right) $ be two Metzler matrices such that: 
 $$\forall\;i\; \forall\;j\quad m_{ij} \leq n_{ij}\quad\mathrm{( \,which \;  we \;denote\;by}\; M \leq N\;\mathrm{or}\; N - M \geq 0 \;\mathrm{)} $$
 
 Denote $\dr{x}(t,\dr{x}_0) = \e^{tM}\dr{x}_0$ and  $\dr{y}(t,\dr{y}_0) = \e^{tN}\dr{y}_0$ then
 $$\forall \; \dr{x}_0\;  \forall \; \dr{y }_0\; \forall \;t\geq 0 \quad 0 \leq  \dr{x}_0 \leq \dr{y}_0 \implic \dr{x}(t,\dr{x}_0)  \leq  \dr{y}(t,\dr{y}_0)$$
 \end{prop}
 \textbf{Proof}. Let $\dr{z}(t) = \dr{y}(t,\dr{y}_0) - \dr{x}(t,\dr{x}_0)$. We have
 $$ \frac{d \dr{z}}{dt} = N\dr{y}(t,\dr{y}_0) - M\dr{x}(t,\dr{x}_0) = N(\dr{y}(t,\dr{y}_0) - \dr{x}(t,\dr{x}_0)) + (N - M)\dr{x}(t,\dr{x}_0)$$
 $$ \frac{d \dr{z}}{dt} =  N\dr{z}(t)  + (N - M)\dr{x}(t,\dr{x}_0)$$
 and from \eqref{nonhomo} and \eqref{solnonhomo} we have 
 \beq \label{somme}
  \dr{z}(t) = \e^{tN} \dr{z}(0) + \int_0^t \e^{(t-\tau)N} u(\tau) d\tau 
  \feq
 with $u(\tau) = (N - M)\dr{x}(\tau,\dr{x}_0)$.
 Since $ \leq  \dr{x}_0 \leq \dr{y}_0 $ we have $\dr{z}(0) \geq 0$ and since  $N$ is Metzler, from Proposition \ref{orthaninvar} we have $\e^{tN}\dr{z}(0) \geq 0$.

 Let us look now to the second term of \eqref{somme}. Since $M$ is Metzler and $\dr{x}_0\geq 0$ we have   $\dr{x}(\tau,\dr{x}_0) \geq 0$ and from the hypothesis $N - M \geq 0$ we have $(N-M) \dr{x}(\tau,\dr{x}_0) = u(\tau) \geq 0$. Since $N$ is  Metzler, again from Proposition \ref{orthaninvar} we have
 $$ \forall  \tau  \leq t\quad \e^{t-\tau)N}u(\tau) \geq 0 $$
 which implies  
 $$  \int_0^t \e^{(t-\tau )N} u(\tau) d\tau   \geq 0 $$
 and achieves the proof of the proposition.$\Box$
%=====================================
 \section{Proof of Proposition \ref{minocircuit}}\label{demprop}
 %=====================================

%============================
 \subsection{Integration along a path.}\label{integrationchemin}
  One considers the system defined on the network \\
\includegraphics[width=1\textwidth]{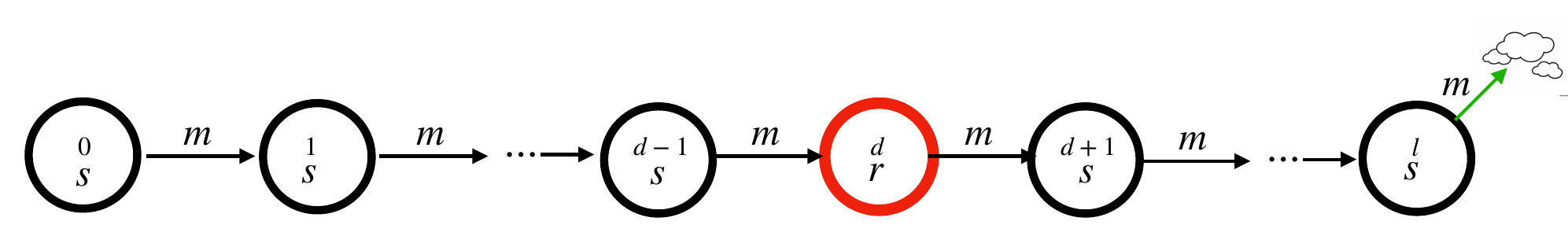}\\
by the equations 
 \beq \label{lemmechain}														
\Sigma(r,s,l,dm,T)\left\{
\begin{array} {ll}
	  j = 0 \cdots d-1& \left\{
	\begin{array}{lcl}
		\displaystyle \frac{dy_0}{dt}& =&  \big(s - m\big) y_{0}\\[8pt]
		\displaystyle \frac{dy_{j}}{dt}& =&  m y_{j-1}+ \big(s - m\big) y_{j}\\[8pt]
		\displaystyle \frac{dy_{d-1}}{dt}& =&  m y_{d-2}+ \big(s-m\big) y_{d-1}
	\end{array} \right.\\[10pt]
 j = d& \left\{
	\begin{array}{lcr}
		\displaystyle \frac{dy_{d}}{dt}\quad & =&  m y_{d-1}+ \big(r-m\big) y_{d}
	\end{array} \right.\\[10pt]
  j = d+1 \cdots l& \left\{
	\begin{array}{lcl}
		\displaystyle \frac{dy_{d+1}}{dt}& =& my_{d} +  \big( s- m\big) y_{d+1}\\[8pt]
		\displaystyle \frac{dy_{j}}{dt}& =&  m y_{j-1}+ \big(s - m\big) y_{j}\\[8pt]
		\displaystyle \frac{dy_{l}}{dt}& =&  m y_{l-1}+ (s-m) y_{l}
	\end{array}
	\right.
\end{array}\right.
\feq
with $r > s$ and the initial condition  $(y_0(t_0) > 0,0\cdots, 0,\cdots,0)$. 
\begin{lemme} \label{lemme}
Consider the function 
\beq \label{cdet1}
C(\theta) =  v_{r,s,d}\left(\frac{\theta}{l-d+1}\right)    \left(w_{r,s}\left(\frac{\theta}{l-d+1}\right)\right)^{l-d}
 \feq
 with
\beq \label{cdet2}
v_{r,s,d}(t) = \int_0^t \e^{-\tau(r-s)}\frac{\tau^{d-1}}{(d-1) !} d\tau, \quad \quad w_{r,s}(t)= \frac{1}{r-s} \left(1-\e^{-t(r-s)}\right).
\feq
Then $\theta \mapsto C(\theta)$ is nondecreasing and one has the following minorization
\beq \label{ineg}
 t \geq   \theta  \implic y_l(t+t_0) \geq C(\theta) m^l \e^{t(\gr{r}-m)}y_0(t_0).
\feq
\end{lemme}

\noindent \textbf{Proof.}

Assume that  $t_0 = 0$.
In the vector form one has
$$\frac{d\dr{y}}{dt} = A\dr{y}$$
with
\begin{small}
\beq
A= \left[\begin{array}{cccccccccccc}
s-m&0&\cdots &\cdots &\cdots&\cdots &\cdots&\cdots&0&0 \\
m&s-m&0 &\cdots&\cdots &\cdots&\cdots&\cdots&\cdots&0 \\
0&m&s-m&0&\cdots&\cdots &\cdots&\cdots&\cdots&0 \\
\cdots&\cdots&\cdots&\cdots&\cdots&\cdots&\cdots&\cdots&\cdots&\cdots\\
0&\cdots&\cdots&m&s-m&\cdots&\cdots&\cdots&\cdots&0\\
0&\cdots&\cdots &\cdots&\red{m}&\red{r-m}&0&\cdots&\cdots&0\\
0&\cdots&\cdots &\cdots&\cdots &m&s-m&0&\cdots&0\\
\cdots&\cdots&\cdots&\cdots&\cdots&\cdots&\cdots&\cdots&\cdots&\cdots\\
\cdots&\cdots&\cdots&\cdots&\cdots&\cdots&\cdots&m&s-m&0 \\
0&\cdots&\cdots &\cdots&\cdots&\cdots &\cdots&\cdots &m&s-m
\end{array}
\right]
\feq
\end{small}
\noindent 
One sees immediately that  $$\dr{x}(t) = \e^{-t(s-m)}\dr{y}(t)$$ is solution of:
$$\frac{d\dr{x}}{dt} = B \dr{x} $$
with
\begin{small}
\beq
B= \left[\begin{array}{cccccccccccc}
0&0&\cdots &\cdots &\cdots&\cdots &\cdots&\cdots&0&0 \\
m&0&0 &\cdots&\cdots &\cdots&\cdots&\cdots&\cdots&0 \\
0&m&0&0&\cdots&\cdots &\cdots&\cdots&\cdots&0 \\
\cdots&\cdots&\cdots&\cdots&\cdots&\cdots&\cdots&\cdots&\cdots&\cdots\\
0&\cdots&\cdots&m&0&\cdots&\cdots&\cdots&\cdots&0\\
0&\cdots&\cdots &\cdots&\red{m}&\red{r}-\red{s}&0&\cdots&\cdots&0\\
0&\cdots&\cdots &\cdots&\cdots &m&0&0&\cdots&0\\
\cdots&\cdots&\cdots&\cdots&\cdots&\cdots&\cdots&\cdots&\cdots&\cdots\\
\cdots&\cdots&\cdots&\cdots&\cdots&\cdots&\cdots&m&0&0 \\
0&\cdots&\cdots &\cdots&\cdots&\cdots &\cdots&\cdots &m&0
\end{array}
\right]
\feq
\end{small}
and the initial condition $\dr{x}(0) = (x_0(0)=y_0(0) > 0,0\cdots, 0,\cdots,0)$.

 \noindent  \textbf{From site $0$ to site $d-1$}.
By successive integrations one has 
 \beq
 x_{d-1}(t) = m^{d-1} \frac{t^{d-1}}{(d-1) !}x_0(0)
 \feq
 \noindent  \textbf{On the site  $d$}. \\
 From Proposition \eqref{propfond}  the integration of the differential equation
$$
\frac{dx_{d}}{dt}  =  m x_{d-1}+ \big(r-s\big) x_{d}\quad \quad x_d(0) = 0
$$ 
gives 
 $$
 x_d(t) = \e^{t(r-s)} \int_0^t \e^{-\tau(r-s)} m x_{d-1}(\tau )d\tau  =  \e^{t(r-s)} \int_0^t \e^{-\tau(r-s)} m^d \frac{\tau^{d-1}}{(d-1) !}x_0d\tau
 $$
 Set 
 \beq
 \boxed{
 v_{r,s,d}(t) = \int_0^t\e^{-\tau(r-s)}\frac{\tau ^{d-1}}{(d-1) !} d\tau}
 \feq
 We have 
 $$
 x_d(t) = \e^{t(r-s)} m^dx_0(0)  v_{r,s,d}(t)
 $$
and, since $v_{r,s,d}(t)$  is an increasing function of $t$ one has for every $\theta >0$ and every integer $p$

 \beq \label{ineg1}
 \boxed{
 t \geq \frac{\theta}{p} \implic x_d(t) \geq  \e^{t(r-s)} m^{d} x_0(0) v_{r,s,d}\left(\frac{\theta}{p} \right)}
 \feq
 where the parameter $p$ will be specified later.\\\\
  \textbf{On the site $d+1$.} \\One has 
 $$
 x_{d+1}(t)  = m \int_0^t x_d(\tau)d\tau 
 $$
 and from \eqref{ineg1} one has
 $$
 t \geq \frac{\theta}{p} \implic  x_{d+1}(t) \geq m\int_{\frac{\theta}{p}}^t  \e^{\tau(r-s)} m^{d} x_0(0) v_{r,s,d}\left(\frac{\theta}{p} \right) d\tau = m^{d+1} x_0(0) v_{r,s,d}\left(\frac{\theta}{p} \right)\int_{\frac{\theta}{p}}^t  \e^{\tau(r-s)}  d\tau 
 $$
 One has
 $$
 \int_{\frac{\theta}{p}}^t  \e^{\tau(r-s)}  d\tau = \frac{1}{r-s}\left( \e^{t(r-s)} - \e^{\frac{\theta}{p}(r-s)}\right) = \frac{1}{r-s} \e^{t(r-s)}\left( 1- \e^{(\frac{\theta}{p}-t)(r-s)}\right)
 $$ 
 And now, since $r-s> 0$, we have 
 \beq
 t \geq 2 \frac{\theta}{p}  \implic   \int_{\frac{\theta}{p}}^t  \e^{\tau(r-s)}  d\tau  \geq \frac{1}{r-s} \e^{t(r-s)}\left( 1- \e^{-\frac{\theta}{p}(r-s)}\right)
 \feq
 If we denote
 \beq
 \boxed{
 w_{r,s}(t) = \frac{1}{r-s}\left( 1- \e^{-t(r-s)}\right)}
 \feq
 it follows
 $$
 t \geq  \frac{\theta}{p} +1\times \frac{\theta}{p} \implic  x_{d+1}(t) \geq m^{d+1} x_0(0) v_{r,s,d}\left(\frac{\theta}{p}\right)  w_{r,s}\left(\frac{\theta}{p}\right) \e^{t(r-s)}
 $$
\beq \label{ineg2}
\boxed{
 t \geq  2 \frac{\theta}{p}  \implic  x_{d+1}(t) \geq   \e^{t(r-s)}   w_{r,s}\left(\frac{\theta}{p}\right)      m^{d+1} x_0(0) v_{r,s,d}\left(\frac{\theta}{p}\right) }
  \feq
\textbf{On the site $d+2$.} \\
One has 
$$
 x_{d+2}(t)  = m \int_0^t x_{d+1}(\tau)d\tau 
 $$
and from \eqref{ineg2}
 $$
 t \geq  2\frac{\theta}{p}  \implic  x_{d+2} \geq w_{r,s}\left(\frac{\theta}{p}\right) m^{d+2} x_0(0) v_{r,s,d}\left(\frac{\theta}{p}\right) \int_{\frac{\theta}{p} +1\times \frac{\theta}{p}}^t \e^{\tau(r-s)}d\tau
 $$
 and like previously 
 $$
 t \geq  3 \frac{\theta}{p} \implic  x_{d+2} \geq w_{r,s}\left(\frac{\theta}{p}\right) m^{d+2} x_0(0) v_{r,s,d}\left(\frac{\theta}{p}\right) w_{r,s}\left(\frac{\theta}{p}\right)\e^{t(r-s)} 
 $$
 \beq
 \boxed{
 t \geq  3 \frac{\theta}{p}  \implic  x_{d+2} \geq \e^{t(r-s)} \left(w_{r,s}\left(\frac{\theta}{p}\right)\right)^2 m^{d+2} x_0(0) v_{r,s,d}\left(\frac{\theta}{p}\right) }
  \feq 
Iterating the process up to $k$ we have on the site $d+k$.
  \beq \label{ineg3}
 \boxed{
 t \geq (k+1)  \frac{\theta}{p}  \implic  x_{d+k} \geq \e^{t(r-s)} \left(w_{r,s}\left(\frac{\theta}{p}\right)\right)^k m^{d+k} x_0(0) v_{r,s,d}\left(\frac{\theta}{p}\right) }
  \feq 
  \textbf{On the site $l =d+k$.}\\
   If we apply \eqref{ineg3} with $k = l-d$ and $p = l-d+1$ we get
  \beq 
  \boxed{
 t \geq  \theta \implic  x_{l} \geq \e^{t(r-s)} \left(w_{r,s}\left(\frac{\theta}{l-d+1}\right)\right)^{l-d} m^{l} x_0(0) v_{r,s,d}\left(\frac{\theta}{l-d+1}\right) }
  \feq 
and if we put $C(\theta) = \left(w_{r,s}\left(\frac{\theta}{l-d+1}\right)\right)^{l-d} v_{r,s,d}\left(\frac{\theta}{l-d+1}\right) $
  \beq 
  \boxed{
 t \geq  \theta \implic  x_{l}(t) \geq C(\theta) \e^{t(r-s)} m^{l} x_0(0)  }
  \feq 
  Now, if we turn back to $y_l(t) = \e^{t(s-m)} x_l(t)$ 
 \beq
t > \theta  \implic  y_{l}(t) > C(\theta)m^l\e^{t(r-m )}y_0(0)
\feq
and, since the matrix $A$ does not depend on $t$, for any  initial condition $y_0(t_0)$ 
 
 \beq
t > \theta  \implic  y_{l}(t+t_0) > C(\theta)m^l\e^{t(r-m )}y_0(t_0)
\feq
which ends the proof of the Lemma.$\Box$

%==============================
 \subsection{Minorization through a path}
 
 Let
$ \Pi = \left\{ \pi_1,\cdots,\pi_i,\cdots,\pi_n \right\}$
be a set of $n $ sites.
In this subsection we consider the system 
\beq \label{static}
\Sigma(r_i,\ell_{ij}m,T)\quad \quad \left \{\frac{dx_i}{dt} = r_i x_i+ m \sum_{j=1}^n \ell_{ij} x_j  \quad i = 1,\cdots,n\right.
\feq
where $r_i$ are constant  and $\ell_{ij} \in \{0, 1\}$ are   associated to the  static network $\CN$ 
$$ \pi_i\to \pi_j  \in \CN \Longleftrightarrow  \ell_{ji} = 1$$
\noindent Consider an arbitrary   simple path $a\Gamma b$ of $\CN$ like the following one \\
\includegraphics[width=0.9\textwidth]{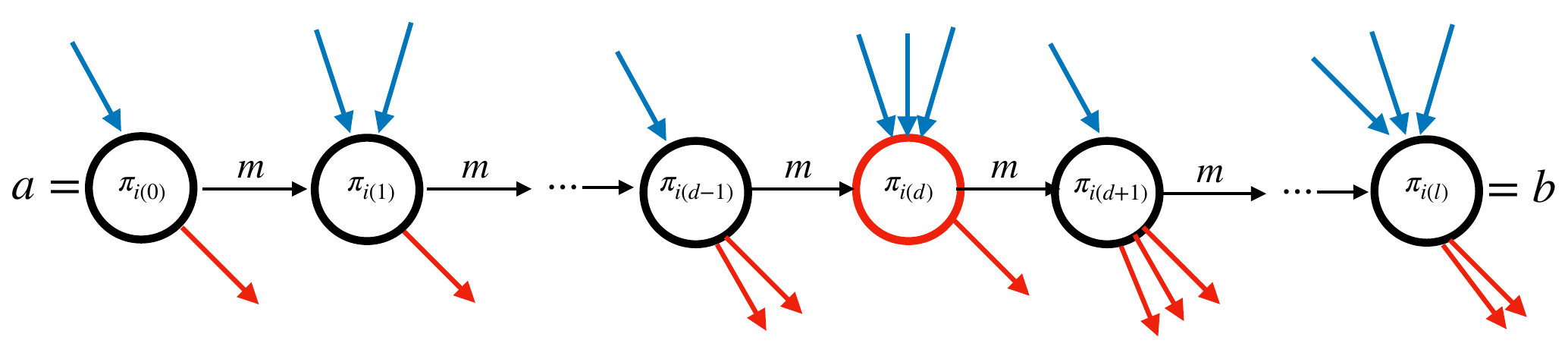}\\
For each site $\pi_{i(j)}$ we have represented in blue the incoming links (from any site of the network) in the site, in black the link that connects to the next site in the path, and finally, in red, the links that leave the site to some other site of the network. Notice that there is no ‘‘black arrow'' leaving the last site and that a red arrow leaving some site might be a link which connect to a site of the chain provided it is not the next one.  \\\\
Our aim in this subsection is to minorize $x_{i(l)}(t) $ given an initial condition such that $x_{i(0)} (t_0) > 0$. \\\\
Compare the following picture to the previous one: \\
\includegraphics[width=1\textwidth]{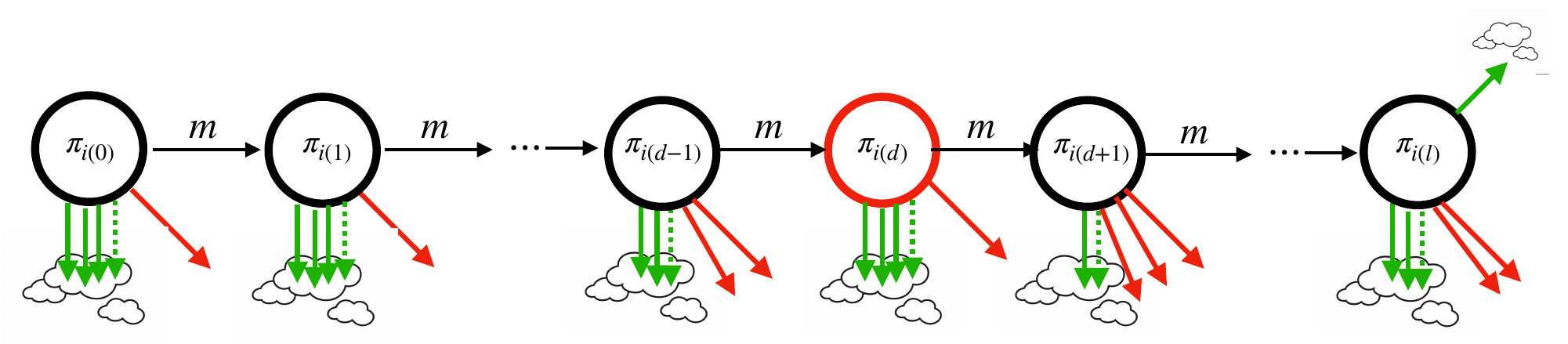}\\
\bito
\item We have cut all the blue links.
\item The number of links leaving each site (red+black) is smaller than $n-1$. We have added links (in green) to the ‘‘clouds'' in a number such that the total links leaving the site  is just the maximum total number $n-1$ of possible links leaving a site. \fit
From this picture we define a new system on $\CN$ in the following manner:
\beq
\Sigma(r_i, a\Gamma b,m,T)\left\{
\begin{array}{rcl}
\displaystyle   \pi_i  \notin a\Gamma b &\implic&\displaystyle \frac{d\xi_i}{dt} = \big(r_i- \alpha- m\big) \xi_i\\[8pt]
\displaystyle i = i(0) &\implic&\displaystyle  \frac{d\xi_{i(0)}}{dt} =\big( r_{i(0)} - \alpha- m \big )\xi_i(0)\\[8pt]
\displaystyle 0 <  j \leq l &\implic &\displaystyle  \frac{d\xi_{i(j)}}{dt} = m \xi_{i(j-1)} + \big( r_{i(j)} - \alpha -m \big)\xi_{i(j)}\\[8pt]
%\displaystyle i = i(l) &\implic& \displaystyle  \frac{d\xi_{i(l)}}{dt} = m \xi_{i(l-1)} + \big( r_{i(l)} - \alpha-m  \big)\xi_{i(j)}
\end{array}
\right.
\feq 
where $\alpha = (n-2)\times m$.
By construction it is evident that 
the system $\Sigma(r_i,a\Gamma b,m,T)$ minorizes  the system $ \Sigma(r_i,\ell_{ij},m,T)$
$$\Sigma(r_i,a\Gamma b,m,T ) \leq  \Sigma(r_i,\ell_{ij},m,T)$$
and that its restriction to $a\Gamma b$, also denoted by $\Sigma(r_i,a\Gamma b,m,T)$, is independent of the $\xi_i(t)$ for $\pi_i \notin a \Gamma b$.

Now define a \textbf{ ‘‘dominant'' site} $i(d)$ as a site such that for every $j = 0\cdot \cdot l$ one has $r_{i(d)}\geq r_{i(j)}$ ; let $\sigma =  \min_{j = 0..l}r_{i(j)} -1$. The ‘‘-1''  in the definition of $\sigma $ is there to ensure that $\sigma < r_{i(d)} $. Put
$$ s = \sigma - \alpha \quad \quad \gr{r} = r_{i(d)}$$
 For $j \not = d$ replace  in  $\Sigma(r_i, a\Gamma b,m,T)$, the term  $r_i-\alpha -m$ by $s- m$ and $r_{i(d)}-m$ by $\gr{r}-m$.
 We obtain 
 \beq \label{lemmechain2}														
\Sigma(r,s,,l,d,m,T)\left\{
\begin{array} {ll}
	  j = 0 \cdots d-1& \left\{
	\begin{array}{lcl}
		\displaystyle \frac{dy_0}{dt}& =&  \big(s - m\big) y_{0}\\[8pt]
		\displaystyle \frac{dy_{j}}{dt}& =&  m y_{j-1}+ \big(s - m\big) y_{j}\\[8pt]
		\displaystyle \frac{dy_{d-1}}{dt}& =&  m y_{d-2}+ \big(s-m\big) y_{d-1}
	\end{array} \right.\\[10pt]
 j = d& \left\{
	\begin{array}{lcr}
		\displaystyle \frac{dy_{d}}{dt}\quad & =&  m y_{d-1}+ \big(\gr{r}-m\big) y_{d}
	\end{array} \right.\\[10pt]
  j = d+1 \cdots l& \left\{
	\begin{array}{lcl}
		\displaystyle \frac{dy_{d+1}}{dt}& =& my_{d} +  \big( s- m\big) y_{d+1}\\[8pt]
		\displaystyle \frac{dy_{j}}{dt}& =&  m y_{j-1}+ \big(s - m\big) y_{j}\\[8pt]
		\displaystyle \frac{dy_{l}}{dt}& =&  m y_{l-1}+ (s-m) y_{l}
	\end{array}
	\right.
\end{array}\right.
\feq
 Again, by construction 
\beq
\Sigma(r,s,l,d ,m,T ) \leq  \Sigma(r_i,\ell_{ij},m,T)
\feq
The system $\Sigma(r,s,l,d,m,T)$ is exactly the system \eqref{lemmechain} of the Lemma \ref{lemme}  and therefore \eqref{ineg} applies 
\beq \label{inegchemin}
 t \geq   \theta  \implic y_l(t+t_0) \geq C(\theta) m^l \e^{t(\gr{r}-m)}y_0(t_0)= C(\theta) m^l \e^{t(r_{i(d)}-\alpha-m)}y_0(t_0)
\feq
with
\beq \label{cdet1}
C(\theta) =  v_{r,s,d}\left(\frac{\theta}{l-d+1}\right)    \left(w_{r,s}\left(\frac{\theta}{l-d+1}\right)\right)^{l-d}
 \feq
 and
\beq \label{cdet2}
v_{r,s,d}(t) = \int_0^t \e^{-\tau(r-s)}\frac{\tau^{d-1}}{(d-1) !} d\tau \quad \quad w_{r,s}(t)= \frac{1}{r-s} \left(1-\e^{-t(r-s)}\right)
\feq
Replacing $r$ and $s$ by their definition one has $r-s = r_{i(d)} -\sigma$ which gives
\beq \label{cdet3}
v_{r,s,d}(t) = \int_0^t \e^{-\tau( r_{i(d)}-\sigma)}\frac{\tau^{d-1}}{(d-1) !} d\tau \quad \quad w_{r,s}(t)= \frac{1}{r_{i(d)}-\sigma} \left(1-\e^{-t( r_{i(d)}-\sigma)}\right)
\feq 
which is independent of $m$.
We have proved the proposition
\begin{prop} \label{minorchemin}
Consider the system $\Sigma(r_i,\ell_{ij},m,T)$
\beq \label{static2}
\Sigma(r_i,\ell_{ij},m,T)\quad \quad \left \{\frac{dx_i}{dt} = r_i x_i+ m \sum_{j=1}^n \ell_{ij} x_j  \quad i = 1,\cdots,n\right.
\feq
defined on some static network $\mathcal{N}$ and a simple path
\beq
a\Gamma b = \{a = \pi_{i(0)} \to\pi_{i(1)} \to \cdots\to \pi_{i(j)}\to \cdots \to \pi_{i(l)} = b\}
\feq  on it. Let 
$r_{i(d)} = \max_{j = 0.. l} r_{i(j)}\quad \quad \sigma =  \min_{j = 0.. l} r_{i(j)}-1$.
Then there exist $\theta$, $C(\theta) > 0$ and $\mu > 0$ such that:
\beq \label{inegchemin2}
 t \geq   \theta  \implic x_l(t+t_0) \geq C(\theta) m^l \e^{t(r_{i(d)}-\mu\times m)}x_0(t_0)
\feq
(take $\mu = n-1$ and  $C(\theta)$ given by \eqref{cdet1} and \eqref{cdet3}).
\end{prop}

\begin{rem}
\label{rem:mu} In our proof we have added red links ‘‘to the clouds'' in number such that the total outgoing links is $n-1$ but it would has been enough in order to use the lemma  to add a number of links such that the number outgoing links is just a constant and have a better minorization. This is why we prefer in the statement of the proposition to be not explicit in the definition of $\mu$.
\end{rem}

%==================================
\subsection{Minorization through a $p$-circuit.}
We come now to the proof of the Proposition \ref{minocircuit}.

Consider the system $\Sigma(r_i^k ,\ell_{ij}^k,m,T)$ on the underlying network \eqref{UN}. Consider the $T-$circuit $\mathcal{C}$ defined by
\beq \label{circuit}
\mathcal{C} =
a_0\Gamma^1 a_1\Gamma^2 a_2 \cdots a_{k-1}\Gamma^k a^k \cdots
a_{p-1}\Gamma^pa_0
\feq
with 
$$ a_{k-1}\Gamma^ka_k = \left\{  a_{k-1} = \pi_{i^k(0)} \to \pi_{i^k(1)}\to \cdots \pi_{i^k(j)} \to \pi_{i^k(l^k)}= a_k\right \}$$
and its growth index $\chi^{\mathcal{C}}$. 

We have to prove that there exist $\theta $ and constants $C >0$, $\mu >0$  (independent of $m$) such that 
\beq \label{minoration}
T > \theta \implic x_{i^1(0)}(T) \geq   C m^{\Lmat } \e^{T\left(\chi^{\mathcal{C}}- \mu \times m\right)} x_{i^1(0)}(0)
\feq
where ${\Lmat }$ is the total length of the circuit.

Let $d^k$ be a dominant site of the path $a_{k-1}\Gamma^k a_k$ and $r^k_{d^k}$ the corresponding dominant rate.\\

\noindent Consider the first simple path $\pi_1\Gamma^1\pi_{i^1(l^1)}$. Let $\theta > 0$ and let 
$C^1  =  C^{\Gamma^1}(\theta)$ where $ C^{\Gamma^1}(\theta)$ is the function of Proposition \ref{minorchemin}  applied to the first path. From \eqref{cdet1} and \eqref{cdet2} we know that $C^1(\theta) > 0$. Let $ T^{*1} \geq \frac{\theta}{t_1}$. From prop \ref{minorchemin} we know that 
\beq \label{circuit1}
T > T^{*1} \implic x_{i^1(l^1)}(Tt_1) > C^1 m^{l^1}\e^{Tt_1\left(r^1_{d^1}-\mu\times m \right)}x_{i^1(0)}(0)
\feq
On the interval $[Tt_1, Tt_2]$ the Proposition \ref{minorchemin} applies again to the simple path $(\pi_{i^1(l^1)}=\pi_{i^2(0)})\Gamma^2 \pi_{i^2(l^2)}$ and thus if we put $C^2 = C^{\Gamma^2}(\theta) $ and  $T^{*2} \geq \frac{\theta}{t_2-t_1}$
\beq \label{circuit2}
 T > T^{*2} \implic x_{i^2(l^2)}(Tt_2) > C^2 m^{l^2}\e^{T(t_2-t_1)\left(r^2_{d^2}-\mu\times m\right)}x_{i^1(l^1)}(Tt_1)
 \feq
Combining the two inequalities \eqref{circuit1} and \eqref{circuit2} we obtain that

$$T > \max\{T^{*1},T^{*2}\} \implic
 x_{i^2(l^2)}(Tt_2) \geq C^1C^2  
 m^{l^1+l^2} \e^{ T\left( (r^1_{d^1}-\mu\times m)t_1+ (r^2_{d^2}-\mu\times m)(t_2-t_1) \right) }x_{i^1(0)}(0).$$

If we iterate this application of Proposition \ref{minorchemin} to the successive simple paths of the circuit up to the last one we get 
$$T> \max\{T^{*1}\cdots T^{*p}\}  \implic x_{i^p(l^p)}(Tt_p) \geq C^1\cdots C^p  m^{\sum_kl^k} \e^{T\left(\sum_{1 =1}^ p r^k_{d^k} (t_k-t_{k-1})-\mu\times m \right)}x_{i^1(0)}(0)$$
Since  $x_{i^p(l^p)}(Tt_p) =x_1(T)$ and by definition of $\chi^{\mathcal{C}}$ one has
$$T> \max\{T^{*1}\cdots T^{*p}\}  \implic x_1(T) \geq C^1\cdots C^p  m^{\sum_kl^k} \e^{T\left(\chi^{\mathcal{C}} -\mu\times m\right)}x_{i^1(0)}(0)$$
which proves the proposition.

 %===================================
 \subsection{Possible extensions}\label{extension}
 The Lemma \ref{lemme} assumes that the path has no loop.
We now give an extension of this lemma in the case of a network with a single simple loop where the dominant site is in the loop.

Consider the system defined by the network:\\
\includegraphics[width=1\textwidth]{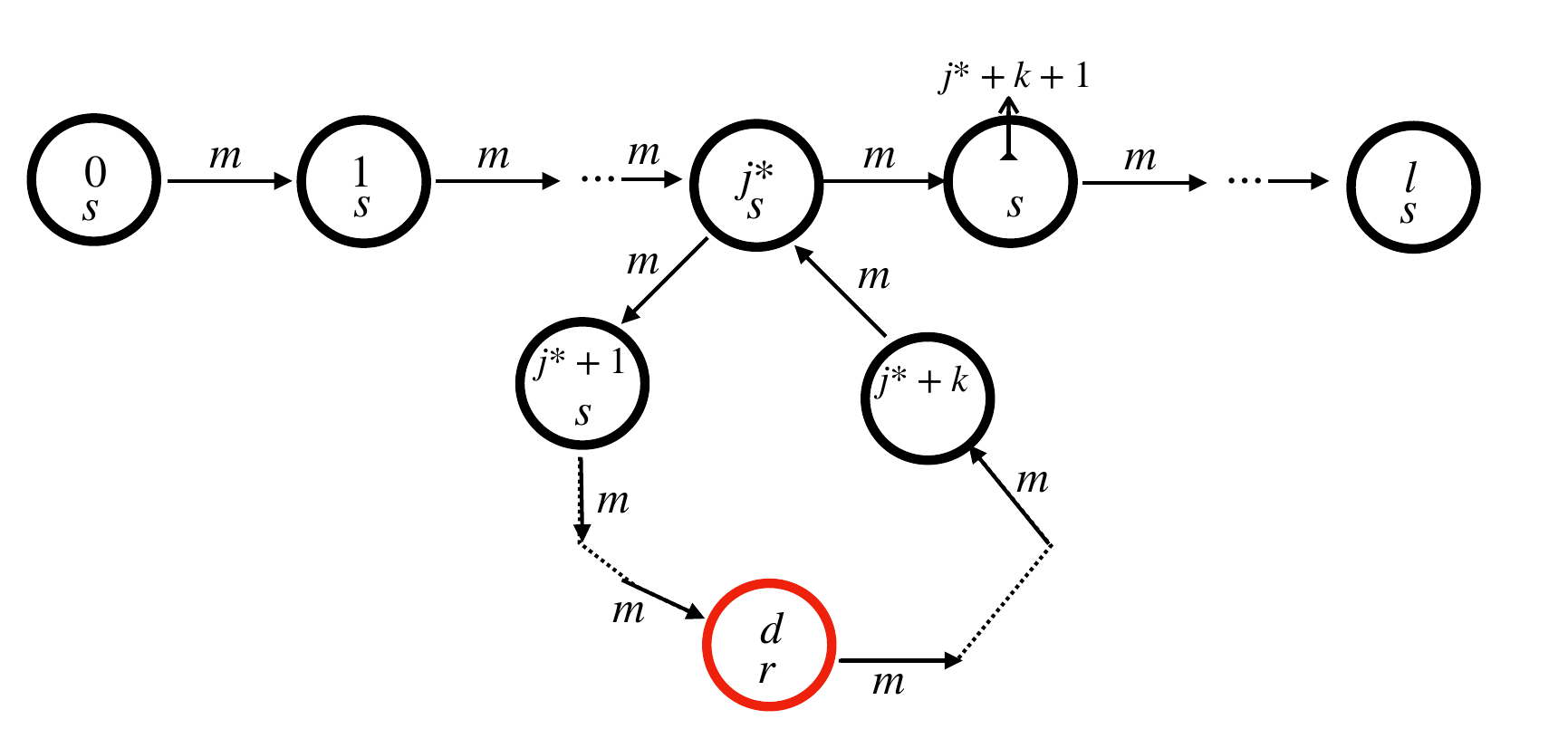}\\
and we leave it to the reader to write down the equations.
\begin{prop} \label{lemme_boucle}
For the above network, we have for all $\theta > 0$ and all $\beta \in (0,1)$, 
\beq
 t \geq   \theta  \implic y_l(t+t_0) \geq C(\beta \theta) C((1 - \beta)\theta) m^{l+1} \e^{t(r-2m)}y_0(t_0)
\feq
\end{prop}

\noindent \textbf{Proof.}

Assume that  $t_0 = 0$. The idea is to cut the loop into two simple paths, to each of which we apply Proposition \ref{minorchemin}. First, consider in the above network the simple path
\[
0 \Gamma d = \{ 0 \to 1 \to \ldots \to j^* \to j^* + 1 \to \ldots \to d \}
\]
Since the maximal number of links leaving a site is $2$ (on site $j^*$), we deduce from Proposition \ref{minorchemin} that for all $\beta \in (0,1)$ and $\theta > 0$, one has
\beq
\label{eq:demiboucle1}
t \geq \beta \theta \implic y_d(t) \geq C( \beta \theta)m^d e^{t (r - 2m)} y_0(0)
\feq
Next, consider the simple path (of  length $l - d + 1$ and maximal index $r$)
\[
d \Gamma l = \{ d \to d+ 1 \to \ldots \to j^* \to j^{*}+k+1 \to j^* + k + l^* = l\}
\]
Once again, the maximal number of links leaving a site is $2$, so that for all $\beta \in (0,1)$, for all $t_0' \geq 0$, for all $\theta > 0$, we have
\beq
\label{eq:demiboucle2}
t' \geq ( 1 - \beta)\theta \implic y_l(t'+t_0') \geq C( ( 1 - \beta) \theta) m^{l - d + 1} e^{t'(r - 2m)} y_d(t_0').
\feq
Taking $t = \beta \theta$ in Equation \eqref{eq:demiboucle1}, and $t_0' = \beta \theta$ in Equation \eqref{eq:demiboucle2} yields
\[
t' \geq ( 1 - \beta)\theta \implic y_l(t'+ \beta \theta) \geq C( ( 1 - \beta) \theta)C( \beta \theta)m^{l+1}e^{(t'+\beta \theta)(r - 2m)}y_0(0)
\]
Noting that $t \geq \theta$ is equivalent to $t' = t - \beta \theta \geq (1 - \beta)\theta$, we finally end up with
\[
t \geq \theta \implic y_l(t) \geq C( ( 1 - \beta) \theta)C( \beta \theta)m^{l+1}e^{t(r - 2m)}y_0(0).
\]

\section*{Acknowledgements}

We thank François Bienvenu and 2 anonymous referees for their valuable comments leading to an improved version of the manuscript.

\section*{Fundings}

\section*{Conflict of interest disclosure}
The authors declare that they comply with the PCI rule of having no financial conflicts of
interest in relation to the content of the article.

\section*{Material and methods}
This article contains only mathematical results with their proofs and interpretations in the field of population dynamics.
%=========================================
{}

\fin